\vsize=9.0in\voffset=1cm
 \looseness=2


\message{fonts,}

\font\tenrm=cmr10
\font\ninerm=cmr9
\font\eightrm=cmr8
\font\teni=cmmi10
\font\ninei=cmmi9
\font\eighti=cmmi8
\font\ninesy=cmsy9
\font\tensy=cmsy10
\font\eightsy=cmsy8
\font\tenbf=cmbx10
\font\ninebf=cmbx9
\font\tentt=cmtt10
\font\ninett=cmtt9

\font\ninesl=cmsl9
\font\eightsl=cmsl8

\font\nineit=cmti9
\font\eightit=cmti8

\skewchar\ninei='177 \skewchar\eighti='177
\skewchar\ninesy='60 \skewchar\eightsy='60

\def\eightpoint{\def\rm{\fam0\eightrm} 
\normalbaselineskip=9pt
\normallineskiplimit=-1pt
\normallineskip=0pt

\textfont0=\eightrm \scriptfont0=\sevenrm \scriptscriptfont0=\fiverm
\textfont1=\ninei \scriptfont1=\seveni \scriptscriptfont1=\fivei
\textfont2=\ninesy \scriptfont2=\sevensy \scriptscriptfont2=\fivesy
\textfont3=\tenex \scriptfont3=\tenex \scriptscriptfont3=\tenex
\textfont\itfam=\eightit  \def\it{\fam\itfam\eightit} 
\textfont\slfam=\eightsl \def\sl{\fam\slfam\eightsl} 

\setbox\strutbox=\hbox{\vrule height6pt depth2pt width0pt}%
\normalbaselines \rm}

\def\ninepoint{\def\rm{\fam0\ninerm} 
\textfont0=\ninerm \scriptfont0=\sevenrm \scriptscriptfont0=\fiverm
\textfont1=\ninei \scriptfont1=\seveni \scriptscriptfont1=\fivei
\textfont2=\ninesy \scriptfont2=\sevensy \scriptscriptfont2=\fivesy
\textfont3=\tenex \scriptfont3=\tenex \scriptscriptfont3=\tenex
\textfont\itfam=\nineit  \def\it{\fam\itfam\nineit} 
\textfont\slfam=\ninesl \def\sl{\fam\slfam\ninesl} 
\textfont\bffam=\ninebf \scriptfont\bffam=\sevenbf
\scriptscriptfont\bffam=\fivebf \def\bf{\fam\bffam\ninebf} 
\textfont\ttfam=\ninett \def\tt{\fam\ttfam\ninett} 

\normalbaselineskip=11pt
\setbox\strutbox=\hbox{\vrule height8pt depth3pt width0pt}%
\let \smc=\sevenrm \let\big=\ninebig \normalbaselines
\parindent=1em
\rm}

\def\tenpoint{\def\rm{\fam0\tenrm} 
\textfont0=\tenrm \scriptfont0=\ninerm \scriptscriptfont0=\fiverm
\textfont1=\teni \scriptfont1=\seveni \scriptscriptfont1=\fivei
\textfont2=\tensy \scriptfont2=\sevensy \scriptscriptfont2=\fivesy
\textfont3=\tenex \scriptfont3=\tenex \scriptscriptfont3=\tenex
\textfont\itfam=\nineit  \def\it{\fam\itfam\nineit} 
\textfont\slfam=\ninesl \def\sl{\fam\slfam\ninesl} 
\textfont\bffam=\ninebf \scriptfont\bffam=\sevenbf
\scriptscriptfont\bffam=\fivebf \def\bf{\fam\bffam\tenbf} 
\textfont\ttfam=\tentt \def\tt{\fam\ttfam\tentt} 

\normalbaselineskip=11pt
\setbox\strutbox=\hbox{\vrule height8pt depth3pt width0pt}%
\let \smc=\sevenrm \let\big=\ninebig \normalbaselines
\parindent=1em
\rm}

\message{fin format jgr}
\magnification=1200
\font\Bbb=msbm10

\def\R{\hbox{\Bbb R}}

\def\Z{\hbox{\Bbb Z}}
\def\pa{\partial}
\def\ep{\varepsilon}
\def\v{\varphi}
\def\b{\backslash}
\vskip 4 mm

\centerline{\bf Estimates for solutions of Burgers type equations and
some applications}
\footnote{}{\ninerm E-mail addresses: henkin@math.jussieu.fr (G.Henkin),
shan@ccas.ru (A.Shananin),}
\noindent
\footnote{}{\ninerm tumanov@uiuc.edu (A.Tumanov)}

\vskip 2 mm
\centerline{{\bf G.M.Henkin}$^1$,\ \ {\bf A.A.Shananin}$^{2,a}$
\footnote{} {\ninerm $^1$ Partially supported by RFBR projet 02-01-00854},\ \
{\bf A.E.Tumanov}$^3$}
\vskip 2 mm

\noindent
$1$ {\ninerm Universit\'e Pierre et Marie Curie,\ \ case 247, F-75252,
Paris, France}

\noindent
$2$ {\ninerm Computing Center, Academy of Science, 117967 Moscow, Russia}

\noindent
$3$ {\ninerm University of Illinois, Urbana, IL61801, USA}

\vskip 4 mm
{\bf Abstract.}
We obtain  precise large time asymptotics for the Cauchy problem for
Burgers type equations satisfying shock profile condition. The proofs are
based on the exact a priori estimates for (local) solutions of these
equations and the result of [7].

\vskip 4 mm
{\bf R\'esum\'e.}
Nous trouvons asymptotiques pr\'ecises en temps grand des solutions
du probl\`eme de Cauchy pour d'\'equations de type de Burgers admettantes
des profils de chocs. Les preuves sont bas\'ees sur les r\'esultats de [7]
et sur l'estimations a priori pr\'ecises des solutions de ces \'equations.

\vskip 4 mm
{\it MSC:} {\ninerm 35J, 35K, 35L}
\vskip 4 mm
{\it Keywords:} {\ninerm Burgers type equations, intermediate asymptotic,
training waves}
\vskip 4 mm
{\it Mot-cl\'es:} {\ninerm Burgers type \'equations, asymptotiques
interm\'ediaires, ondes progressives.}

\vskip 4 mm
{\bf 1. Introduction.}
The Burgers type equations have been introduced for studing   different
models of fluids ([1],[3],[4],[10]). The difference-differential analogues of these
 equations have been proposed in some models of economic development
([5],[6]).

One of the most useful versions of the Burgers type equations is the
following

\noindent
([4],[11],[13])
$${\pa f\over \pa t}+\v(f){\pa f\over \pa x}=\ep {\pa^2f\over \pa x^2},
\eqno(1.1)$$
where $\ep>0$, $(x,t)\in\Omega\subset\R^2$.

One of the most interesting difference-differential analogues of equation
(1.1) is the following ([5],[6])
$${\pa F\over \pa t}+\v(F){{F(x,t)-F(x-\ep,t)}\over \ep}=0,\eqno(1.2)$$
where $\ep>0$, $(x,t)\in\Omega\subset\R^2$.

The interesting and difficult problems, related with equations (1.1), (1.2),
are the following.

\vskip 2 mm
{\bf Problem I} ([4],[11]). Find asymptotic $(t\to\infty)$ of the
solution $f(x,t)$, $x\in\R$, $t\ge t_0$, of the equation (1.1) with initial
condition:
$$\alpha\le f(x,t_0)\le\beta,\ \ \int_{-\infty}^0(f(x,t_0)-\alpha)dx+
\int_0^{\infty}(\beta-f(x,t_0))dx<\infty.\eqno(1.3)$$

\vskip 2 mm
{\bf Problem II} ([6]). Find asymptotic $(t\to\infty)$ of the solution
$F(n,t)$, $n\in\Z$, $t\ge t_0$, of the equation (1.2) with $\ep=1$ and
initial condition
$$\alpha\le F(n,t_0)\le\beta,\ \ \sum_{-\infty}^0(F(n,t_0)-\alpha)+
\sum_0^{\infty}(\beta-F(n,t_0))<\infty.\eqno(1.4)$$

See [7], [13] for a review of several recent results on these problems.

In this paper we present a complete solution of these problems for the
 special case of equations, satisfying the shock profile condition.
The detailed study of this special case is highly important for
solving these problems (see [6], [7]).

\vskip 2 mm
{\bf Definition.}
The equation (1.1) (correspondingly (1.2)) satisfies $(\alpha,\beta)$-shock
profile condition, if there exist wave-train solutions of this equation of
the form $f={\tilde f}(x-Ct)$ (corr. $F={\tilde F}(x-Ct)$) such that
${\tilde f}(x)\to\beta$, $x\to +\infty$, ${\tilde f}(x)\to\alpha$,
$x\to -\infty$
(correspondingly
${\tilde F}(x)\to\beta$, $x\to +\infty$, ${\tilde F}(x)\to\alpha$,
$x\to -\infty$).

From the results of [4],[12] it follows that equation (1.1)
 with positive $\v$ satisfies (0,1)-shock profile condition iff
$${1\over u}\int_0^u\v(y)dy>C=\int_0^1\v(y)dy,\ \ \forall\ u\in (0,1).
\eqno(1.5)$$
From the results of [5],[2] it follows that equation (1.2) with positive
$\v$ satisfies (0,1)-shock profile condition iff
$${1\over u}\int_0^u{dy\over \v(y)}< {1\over C}=\int_0^1{dy\over \v(y)},\ \
\forall\ u\in (0,1).\eqno(1.6)$$

Let further $\v$ be a positive piecewise twice continuously differential
function on the interval [0,1].

\vskip 2 mm
{\bf Theorem 1.}

i) {\it Let equation} (1.1) {\it satisfy} (0,1)-{\it shock profile condition}
 (1.5);
$\v^{\prime}(0)\ne 0$ {\it if} $\v(0)=C$; $\v^{\prime}(1)\ne 0$ {\it if}
 $\v(1)=C$.
{\it Let} $f(x,t)$ {\it be a solution of} (1.1) {\it with initial condition}
 (1.3), {\it where}
$\alpha=0$, $\beta=1$. {\it Then there exist constants} $\gamma_0$ {\it and}
 $d_0$ {\it such that}
$$\sup_{x\in\R}|f(x,t)-{\tilde f}(x-Ct+\ep\gamma_0\ln\,t+d_0)|\to 0,\
\ t\to\infty,\eqno(1.7)$$
{\it where} ${\tilde f}(x-Ct)$ {\it is a wave-train solution of} (1.1),
$$\gamma_0=\left\{\matrix{
0,\ \ &{\it if}\ \ \v(0)>C>\v(1),\hfill\cr
{1\over \v^{\prime}(1)},\ \ &{\it if}\ \ \v(0)>C=\v(1),\cr
-{1\over \v^{\prime}(0)},\ \ &{\it if}\ \ \v(0)=C>\v(1),\cr
{1\over \v^{\prime}(1)}-{1\over \v^{\prime}(0)},\ \ &{\it if}\ \
\v(0)=C=\v(1).\cr}\right.$$

ii) {\it Let equation} (1.2) {\it satisfy} (0,1)-{\it shock profile condition}
 (1.6);
$\v^{\prime}(0)\ne 0$ {\it if} $\v(0)=C$; $\v^{\prime}(1)\ne 0$ {\it if}
 $\v(1)=C$.
{\it Let} $F(n,t)$ {\it be a solution of} (1.2) {\it with initial condition}
 (1.4), {\it where}
$\alpha=0$, $\beta=1$ {\it and} $\Delta F(n,t_0)\ge 0$
{\it Then there exist constants} $\Gamma_0$ {\it and}
 $D_0$ {\it such that}
$$\sup_{n\in\Z}|F(n,t)-{\tilde F}(n-Ct+\Gamma_0\ln\,t+D_0)|\to 0,\
\ t\to\infty,\eqno(1.8)$$
{\it where}   ${\tilde F}(x-Ct)$ {\it is a wave-train solution of} (1.2),
$\Delta F(n,t)\buildrel \rm def \over = F(n,t)-F(n-1,t)$
$$\Gamma_0=\left\{\matrix{
0,\ \ &{\it if}\ \ \v(0)>C>\v(1),\hfill\cr
{C\over 2\v^{\prime}(1)},\ \ &{\it if}\ \ \v(0)>C=\v(1),\cr
-{C\over 2\v^{\prime}(0)},\ \ &{\it if}\ \ \v(0)=C>\v(1),\cr
{C\over 2}\bigl({1\over \v^{\prime}(1)}-{1\over \v^{\prime}(0)}\bigr),\ \
&{\it if}\ \ \v(0)=C=\v(1).\cr}\right.$$

\vskip 2 mm
{\bf Remarks.}

1. In the case $\v(0)>C>\v(1)$ the statement i) of Theorem 1 is the main
result of [9] and the statement ii) of Theorem 1 is the main result of [5].

2. For the other cases when $\v(0)=C$ or $\v(1)=C$ or $\v(0)=\v(1)=C$ in
the previous work [7] it was already obtained the existence of the
shift-functions $\gamma(t)=O(\ln\,t)$ and $\Gamma(t,\{x\})=O(\ln\,t)$ with the
 properties
$$\eqalign{
&\sup_x|f(x,t)-{\tilde f}(x-Ct+\ep\gamma(t))|\to 0\ \ {\rm and}\cr
&\sup_x|F(x,t)-{\tilde F}(x-Ct+\ep\Gamma(t,\{x\}))|\to 0,\ \ t\to\infty,\cr}$$
where $f,F$ - solutions of (1.1), (1.2) under conditions (1.3), (1.4),
$\{x\}$ is the fractional part of $x\in\R$.

3. It is interesting to compare the statements i), ii) of Theorem 1 with
the $L^1$- stability results presented in the paper of D.Serre [13].
Results of [13] give in particular the following.

Let $f(x,t)$ and $F(n,t)$ be  solutions of equations (1.1) and (1.2)
correspondingly with such initial conditions that
$$\int_{-\infty}^{\infty}|f(x,0)-{\tilde f}(x)|dx<\infty,\ \
\sum_{-\infty}^{\infty}|F(n,0)-{\tilde F}(n)|<\infty,$$
where ${\tilde f}(x-Ct)$ and ${\tilde F}(n-Ct)$ are wave-trains solutions
of (1.1), (1.2). Then
$$\int_{-\infty}^{\infty}|f(x,t)-{\tilde f}(x-Ct+d_0)|dx\to 0,\ \
\sum_{-\infty}^{\infty}|F(n,t)-{\tilde F}(n-Ct+D_0)|\to 0,\ \ t\to\infty,$$
where constants $d_0$ and $D_0$ are being calculated from equations
$$\int_{-\infty}^{\infty}(f(x,0)-{\tilde f}(x+d_0))dx=0\ \ {\rm and}\ \
\sum_{-\infty}^{\infty}\int_{{\tilde F}(n+D_0)}^{F(n,0)}
{dy\over \v(y)}=0.$$

The proof of Theorem 1 is based on the results of [7] and the following
crucial a priori estimates of (local) solutions for (1.1) and (1.2).

Without loss of generality we will put further parameter $\ep$ equal to 1.
Otherwise, we make substitution: $t\to {t\over \ep}$, $x\to {x\over \ep}$.

\vskip 2 mm
{\bf Theorem 2.}

{\it Let in} (1.1), (1.2) {\it parameter} $\ep=1$. {\it Let} $C=\v(0)>0$,
$\gamma_0>|\v^{\prime}(0)|$,
$\bar x\buildrel \rm def \over = {{x-Ct}\over \sqrt{Ct}}$,
$$\Omega_{\sigma}=\{(x,t)\ \ :\ \ a_1<\bar x<a_2+\sigma\sqrt{Ct}\},\ \
0<a_1<a_2<\infty,\ \ \sigma\ge 0.$$

i) {\it If function} $f(x,t)$ {\it defined in the domain} $\Omega_0$
{\it satisfies equation} (1.1) {\it and}
$$|f(x,t)|\le {\gamma\over \sqrt{Ct}},\ \ (x,t)\in\Omega_0,\ \ t\ge t_0,
\eqno(1.9)$$
{\it then the following  estimate holds}
$$\big|{\pa f\over \pa x}(x,t)\big|\le {b\gamma\over Ct},\ \ (x,t)\in\Omega_0,
\ \ t\ge t_0,\eqno(1.10)$$
{\it where}
$$b={b_0\over C}\bigl(\gamma\gamma_0+{1\over \delta}\bigr)
\biggl(1+\ln_+{{\gamma\gamma_0+1/\delta}\over \sqrt{C}}\biggr),$$
$d=\min\,(\bar x -a_1,a_2-\bar x,a_2/2)$, $\delta=\min\,(1,d)$, $b_0$ {\it is
absolute constant}.

ii) {\it If function} $F(x,t)$ {\it defined in the domain} $\Omega_{\sigma}$,
$\sigma>0$,
{\it satisfies equation} (1.2), $\Delta F(x,t)\buildrel \rm def \over =
F(x,t)-F(x-1,t)\ge 0$, $t\ge t_0$ {\it and}
$$|F(x,t)|\le {\Gamma\cdot \bar x\over \sqrt{Ct}},\ \ {\it where}\ \
\bar x\in (a_1,a_2+\sigma\sqrt{Ct}),\ \ t\ge t_0,\eqno(1.11)$$
{\it then the following  estimate holds}
$$0\le\Delta F(x,t)\le {B\Gamma\cdot \bar x\over Ct},\ \ {\it where}\ \
\bar x\in (a_1,a_2+\sigma_0\sqrt{Ct}), \sigma>\sigma_0,\ t\ge
t_0\ge a_1^2,\eqno(1.12)$$
$$B=B_0\bigl[{\sqrt{1+\sigma}\over \sqrt{\sigma-\sigma_0}}+
{\gamma_0\Gamma\over C}+
{1\over d}+{\gamma_0\Gamma\cdot a_1 \over C}\bigr],\ d=\bar x -a_1,$$
$B_0$ {\it is absolute constant}.

$ii)^{\prime}$ {\it If function} $F(x,t)$ {\it defined in the domain}
$\Omega_0$,
{\it satisfies equation} (1.2), $\v^{\prime}(0)\ge 0$,
 $\Delta F(x,t)\ge 0$, $t\ge t_0$ {\it and}
$$0\le |F(x,t)|\le {\Gamma\over \sqrt{Ct}},\ \ (x,t)\in\Omega_0,\ \ t\ge t_0,
\eqno(1.11)^{\prime}$$
$$\Delta F(x,t_0)\buildrel \rm def \over = F(x,t_0)-F(x-1,t_0)\ge 0,$$
{\it then the following  estimate holds}
$$0\le\Delta F(x,t)\le {B\Gamma\over Ct},\ \ (x,t)\in\Omega_0,\ \ t\ge t_0,
\eqno(1.12)^{\prime}$$
{\it where}
$$B=B_0\bigl[a_2+
\bigl({1\over d}+{\gamma_0\Gamma\over C}\bigr)(1+\ln\,(1+a_2))
\bigr],\ d=\min\,(\bar x -a_1,a_2-\bar x),$$
$B_0$ {\it is absolute constant}.

\vskip 2 mm
{\bf Remarks.}

1. Theorem $2ii)^{\prime}$ in the weak form (condition $0\le F(x,t)\le
O(1/\sqrt{t})$ for $\bar x\in [a_1,a_2]$ implies the estimate
$0\le\Delta F(x,t)\le O(1/t)$ for $\bar x\in [{\tilde a}_1,{\tilde a}_2]
\subset (a_1,a_2))$
was formulated in [7] (with the reference to the present paper) and was
essentially used in [7].

2. Theorem 2ii) is used for the proof of Theorem 1ii) of this paper.
Theorem 2i) is needed for the proof of Theorem 1i).

3. Theorem 2 can be applied to the problems I,II, because the necessary
conditions (1.9), (1.11) are always satisfied due to [14], [6].

4. Theorem 2 can be applied also to the study of Problems I,II in  other
cases. For example, in the important case $\alpha=\beta=0$ the necessary
conditions (1.9), (1.11) are valid globally:
$|f(x,t)|=O(1/\sqrt{t})$, $|F(x,t)|=O(1/\sqrt{t})$, $x\in\R$, $t>0$
(see [8],[14],[6]).

Theorem 1ii) is proved in Section 2. The proofs of
Theorem 2ii) and sketch of the proof of Theorem $2ii)^{\prime}$
 are given  in Section 3.
Theorem 1i) and Theorem 2i) will be proved in the another paper.

\vskip 4 mm
{\bf 2. Asymptotics for solutions of Burgers type equations with shock
profile conditions.}
\vskip 2 mm
The detailed proof of Theorem 1ii) will be given below only in the
principal case: $\alpha=0$, $\beta=1$, $\ep=1$, $\v(0)>C=\v(1)$, $x=n\in\Z$.
Other cases can be proved by very similar arguments.

Let $F(n,t)$, $n\in\Z$, $t\in\R_+$, be a solution of the equation
$${dF(n,t)\over dt}=\v(F(n,t))(F(n-1,t)-F(n,t)),\eqno(2.1)$$
under initial conditions: $F(n-1,t_0)\le F(n,t_0)$, $n\in\Z$,
$$\sum_{-\infty}^0F(n,t_0)+\sum_0^{\infty}(1-F(n,t_0)<\infty.\eqno(2.2)$$
By the shock profile condition there exists a wave-train solution
${\tilde F}(n-Ct)$ for (2.1) with overfall (0,1).

Let $\Phi(F)=\int_F^1dy/\v(y)$. Let $d_A(t)$, $A>0$, be such function that

$$\eqalign{
&\sum_{k=-\infty}^{[Ct+A\sqrt{t}]}(\Phi(F(k,t)-\Phi({\tilde F}(k-Ct+d_A(t)))+
(Ct+A\sqrt{t}-[Ct+A\sqrt{t}])\times\cr
&(\Phi(F([Ct+A\sqrt{t}]+1,t))-\Phi({\tilde F}([Ct+A\sqrt{t}]+1-Ct+d_A(t)))=0.
\cr}\eqno(2.3)$$
By Theorem 1 from [7] for any $A>2\sqrt{C}$ we have
$${\Gamma_-\over t}<d_A^{\prime}(t)\buildrel \rm def \over = {d\over dt}
d_A(t)\le {\Gamma_+\over t},\eqno(2.4)$$
where $0<\Gamma_-\le\Gamma_+<\infty$, $t>t_0>0$ and
$$\sup_n|F(n,t)-{\tilde F}(n-Ct+d_A(t))|\to 0,\ \ t\to\infty.\eqno(2.5)$$

To prove Theorem 1ii) we use statement (2.5) and the following crucial
improvement of the statement (2.4).

\vskip 2 mm
{\bf Proposition 1.}
{\it Let} $A>2\sqrt{C}$. {\it Then the shift-function, defined by} (2.3),
{\it has the following asymptotic behavior}
$$d_A(t)={1\over 2}{C\over \v^{\prime}(1)}\ln\,t + const +o(1),\ \ t\to\infty.
\eqno(2.6)$$

The proof of Proposition 1 is based on the appropriate comparison of
statements for Burgers type equations and on Theorem 2ii) proved in Section 3.
Besides known comparison results [5],[7] we need also the following
new one.

\vskip 2 mm
{\bf Lemma 1.}
{\it Let}
$$\psi(z)={C\over \v^{\prime}(1)}\exp\,\bigl(-{z^2\over 2}\bigr)
\biggl(\int_{-\infty}^{z/2}\exp\,(-2y^2)dy\biggr)^{-1}.$$
{\it For any solution} $F(n,t)$ {\it of the Cauchy problem} (2.1),(2.2)
{\it and for any} $0<\delta_0<\delta<1$ {\it and} $A>2\sqrt{C}$ {\it there
exist} $t_0>0$, $T>0$, {\it such that}
$$F(n,t-T)>1-{1\over \sqrt{t}}\psi\bigl({{n-Ct-2\sqrt{Ct}-\delta\sqrt{Ct}}
\over \sqrt{Ct}}\bigr),\eqno(2.7)$$
{\it if} $Ct+2\sqrt{Ct}+(\delta-\delta_0)\sqrt{Ct}<n<Ct+A\sqrt{Ct}$, $t>t_0$.

\vskip 2 mm
{\bf Remark.}
The function $u(\xi,t)=1-{1\over \sqrt{t}}\psi\bigl({\xi\over \sqrt{t}}
\bigr)$ is one of the most important (in fluid mechanics) solutions of
the classical Burgers equation: ${\pa u\over \pa t}+u{\pa u\over \pa \xi}=
{1\over 2}{\pa^2u\over \pa \xi^2}$ (see [10]).

For the proving Lemma 1 we need additional lemmas about subsolutions for the
equation (2.1) and about patching of these subsolutions.

The next lemma shows that the function $1-{1\over \sqrt{t}}
\psi\bigl({{x-Ct}\over \sqrt{Ct}}\bigr)$, being a solution of classical
Burgers equation, is also the subsolution for the equation (2.1) in the
domains
$$\{(x,t)\ \ :\ \ B<{{x-Ct}\over \sqrt{Ct}}<A,\ \ t>t_0\},\ \
t_0=t_0(A,B).$$
This subsolution will be called asymptotic subsolution.

\vskip 2 mm
{\bf Lemma 2.}
{\it For any} $B<A$ {\it and increasing function} $D(t)=O(\sqrt{t})$ there
exists $t_0>0$ {\it such that for} $t\ge t_0$ {\it and}
$x\in (Ct+B\sqrt{Ct},Ct+A\sqrt{Ct})$ {\it the function}
${\hat F}(x,t)=1-{1\over \sqrt{t}}\psi\bigl({{x-Ct-D(t)}\over \sqrt{Ct}}
\bigr)$ {\it satisfies inequality}
$${\pa {\hat F}\over \pa t}(x,t)\le\v({\hat F}(x,t))({\hat F}(x-1,t)-
{\hat F}(x,t)).\eqno(2.8)$$

\vskip 2 mm
{\bf Remark.}

For the proof of Lemma 1 we will use Lemma 2 in the domain

\noindent
$\{(x,t)\ \ :\ \ 2-\delta<{{x-Ct}\over \sqrt{Ct}}<A\}$ for
$D(t)=(2+\delta_0)\sqrt{Ct}$, $\delta_0<\delta<1$.

In other domains
$\{(x,t)\ \ :\ \ 1<{{x-Ct}\over \sqrt{Ct}}\le 2-\delta\}$ and
$\{(x,t)\ \ :\ \ {{x-Ct}\over \sqrt{Ct}}\le 1\}$ we will need other
subsolutions for (2.1): so called diffusion subsolution
${\hat F}(x,t)=\v^{(-1)}\bigl({{x-2\sqrt{Ct}}\over t}\bigr)$ and
wave-train subsolution ${\tilde F}_{\sigma}(x-C_{\sigma}t)$ with overfall
$[-\sigma,1]$, $\sigma>0$ (see the properties of these subsolutions in
[5],[6]).

\vskip 2 mm
{\bf Proof of Lemma 2.}
We will use the equality
$${\pa {\hat F}(x,t)\over \pa t}={1\over 2t^{3/2}}{\hat \psi}\bigl({{x-Ct}
\over 2\sqrt{t}}\bigr)-{1\over \sqrt{t}}{d{\hat\psi}\bigl({{x-Ct}\over
2\sqrt{t}}\bigr)\over d\bar x}\cdot \bigl(-{x\over 4t^{3/2}}-
{C\over 4\sqrt{t}}\bigr),$$
where
$$\hat\psi(\bar x)={C\over \v^{\prime}(1)}\exp\,(-{2\over C}{\bar x}^2)
\biggl(\int_{-\infty}^{\bar x}\exp\,(-{2\over C}y^2)dy\biggr)^{-1},\ \
\bar x={{x-Ct}\over 2\sqrt{t}}.$$

Let us fixe $\beta>0$. Then for $\bar x={{x-Ct}\over 2\sqrt{t}}\ge -\beta$
and  $t\to +\infty$ we have
$$\eqalign{
&\v({\hat F}(x,t))=C-{\v^{\prime}(1)\over \sqrt{t}}
{\hat \psi}\bigl({{x-Ct}\over 2\sqrt{t}}\bigr)+O\bigl({\psi^2(-\beta)\over
t}\bigr),\cr
&{\hat F}(x-1,t)-{\hat F}(x,t)=-{\pa {\hat F}(x,t)\over \pa x}+{1\over 2}
{\pa^2{\hat F}\over \pa x^2}(x,t)+\ldots =\cr
&2\bigl({1\over 2\sqrt{t}}\bigr)^2
{d{\hat\psi}\bigl({{x-Ct}\over 2\sqrt{t}}\bigr)\over d\bar x}-\bigl({1\over
2\sqrt{t}}\bigr)^3
{d^2{\hat\psi}\bigl({{x-Ct}\over 2\sqrt{t}}\bigr)\over d{\bar x}^2} +
O(1/t^2).\cr}$$
Hence, for $\bar x\ge -\beta$ we obtain
$$\eqalign{
&{\pa {\hat F}(x,t)\over \pa t}-\v({\hat F}(x,t))({\hat F}(x-1,t)-{\hat F}
(x,t))=\cr
&{1\over 2t^{3/2}}{\hat\psi}(\bar x)+{1\over t^2}{d{\hat\psi}(\bar x)\over
d\bar x}\bigl({{2Ct+2{\bar x}\sqrt{t}}\over 4}\bigr)-\cr
&(C-{\v^{\prime}(1)\over \sqrt{t}}{\hat\psi}(\bar x))
\bigl({1\over 2t}{d{\hat\psi}(\bar x)\over d{\bar x}}-{1\over 8t^{3/2}}
{d^2{\hat\psi}(\bar x)\over d {\bar x}^2}\bigr)+O(1/t^2)=\cr
&{1\over 2t^{3/2}}\biggl({\hat\psi}(\bar x)+{\bar x}{d{\hat\psi}(\bar x)\over
d{\bar x}}+\v^{\prime}(1){\hat\psi}(\bar x){d{\hat\psi}(\bar x)\over
d{\bar x}}+{C\over 4}{d^2{\hat\psi}(\bar x)\over d{\bar x}^2}\biggr)+O(1/t^2).\cr}
\eqno(2.9)$$

By direct differentiation with respect to $\bar x$ we obtain
$${d^2{\hat\psi}\over d{\bar x}^2}+\bigl({4\over C}{\bar x}+
{2\v^{\prime}(1)\over C}{\hat\psi}(\bar x)\bigr){d{\hat\psi}(\bar x)\over
d{\bar x}}+{4\over C}{\hat\psi}(\bar x)=0.$$
Hence,
$$\eqalign{
&{\hat\psi}(\bar x)+{\bar x}{d{\hat\psi}\over d{\bar x}}+\v^{\prime}(1)
{\hat\psi}(\bar x){d{\hat\psi}(\bar x)\over d{\bar x}}+{C\over 4}
{d^2{\hat\psi}(\bar x)\over d{\bar x}^2}=\cr
&{\hat\psi}(\bar x)+{\bar x}{d{\hat\psi}\over d{\bar x}}+\v^{\prime}(1)
{\hat\psi}(\bar x){d{\hat\psi}(\bar x)\over d{\bar x}}-
{\bar x}{d{\hat\psi}\over d{\bar x}}-{\v^{\prime}(1)\over 2}
{\hat\psi}(\bar x){d{\hat\psi}(\bar x)\over d{\bar x}}-{\hat\psi}(\bar x)=\cr
&{\v^{\prime}(1)\over 2}{\hat\psi}(\bar x){d{\hat\psi}(\bar x)\over
d{\bar x}}.\cr}$$

Let us check the inequality
$${d{\hat\psi}(\bar x)\over d{\bar x}}<0,\ \ \forall\bar x\in\R.\eqno(2.10)$$
By direct differentiation we have
$${d{\hat\psi}(\bar x)\over d{\bar x}}=-{4\over C}{\bar x}
{\hat\psi}(\bar x)- {\v^{\prime}(1)\over C}{\hat\psi}^2(\bar x).$$
This implies the  equality
$$\eqalign{
&{d{\hat\psi}(\bar x)\over d{\bar x}}=\cr
&-{\hat\psi}(\bar x)
\biggl(\int_{-\infty}^{\bar x}\exp\,(-{2\over C}y^2)dy\biggr)^{-1}
\biggl(\int_{-\infty}^{\bar x}\exp\,(-{2\over C}y^2)dy+
{C\exp\,(-2{\bar x}^2/C)\over 4{\bar x}}\biggr){4\bar x\over C}.\cr}$$
Hence, (2.8) is equivalent to the inequality
$${4\over C}{\bar x}\int_{-\infty}^{\bar x}\exp\,(-{2\over C}y^2)dy+
\exp\,(-2{\bar x}^2/C)>0.$$
For $\bar x\ge 0$ this inequality is obvious. For $\bar x<0$
this inequality follows from the relations
$$\eqalign{
&\lim_{{\bar x}\to -\infty}
\biggl(\int_{-\infty}^{\bar x}\exp\,(-{2\over C}y^2)dy+
{C\exp\,(-2{\bar x}^2/C)\over 4{\bar x}}\biggr)=0\ \ {\rm and}\cr
&{d\over d{\bar x}}
\biggl(\int_{-\infty}^{\bar x}\exp\,(-{2\over C}y^2)dy+
{C\exp\,(-2{\bar x}^2/C)\over 4{\bar x}}\biggr)=-
{C\exp\,(-2{\bar x}^2/C)\over 4{\bar x}^2}<0.\cr}$$.

From (2.9), (2.10) it follows that there exists $\sigma>0$ such that
$$\sup\,\bigl\{{\hat\psi}(\bar x){d{\hat\psi}(\bar x)\over d{\bar x}}\ \
\bigg| \ \ -\beta\le \bar x\le\alpha\bigr\}<-\sigma.$$
Hence, for $x\in [Ct-\beta\sqrt{t},Ct+\alpha\sqrt{t}]$ we obtain the
estimate:
$${\pa {\hat F}(x,t)\over \pa t}-\v({\hat F}(x,t)({\hat F}(x-1,t)-
{\hat F}(x,t))\le -{\v^{\prime}(1)\over 4t^{3/2}}\sigma+O(1/t^2).$$

It means that there exists $t_0>0$ such that for $t\ge t_0$ and
$x\in [Ct-\beta\sqrt{t},Ct+\alpha\sqrt{t}]$ the inequality (2.8) is valid if
${\hat F}(x,t)=1-{1\over \sqrt{t}}\psi\bigl({{x-Ct}\over \sqrt{Ct}}\bigr)$.
This inequality is also valid if
${\hat F}(x,t)=1-{1\over \sqrt{t}}\psi\bigl({{x-Ct-D(t)}\over \sqrt{Ct}}
\bigr)$
for   $x\in [Ct+D(t)-\beta\sqrt{t},Ct+D(t)+\alpha\sqrt{t}]$ because
$t\mapsto D(t)$ is increasing function. The next lemma gives conditions
for patching diffusion subsolutions and asymptotic subsolutions.

\vskip 2 mm
{\bf Lemma 3.}
{\it For any} $\delta\in (0,1)$ {\it and constant} $\Gamma>0$ {\it there
exists} $t_0>0$ {\it such that for} $t\ge t_0$ {\it and}
$n\in [Ct+(2-\delta)\sqrt{Ct}-\Gamma,Ct+(2-\delta)\sqrt{Ct}+\Gamma]$ {\it the
following inequality is valid}
$$\v^{(-1)}\bigl({{n-2\sqrt{Ct}}\over t}\bigr)>1-{1\over \sqrt{t}}
\psi\bigl({{n-Ct-2\sqrt{Ct}}\over \sqrt{Ct}}\bigr).\eqno(2.11)$$

\vskip 2 mm
{\bf Proof of Lemma 3.}
We have equalities
$$\eqalign{
&\lim_{{\bar x}\to -\infty}{1\over {\bar x}}\exp\,(-2{\bar x}^2)
\biggl(\int_{-\infty}^{\bar x}\exp\,(-2y^2)dy\biggr)^{-1}=-4\ \ {\rm and}\cr
&\lim_{{\bar x}\to -0}{1\over {\bar x}}\exp\,(-2{\bar x}^2)
\biggl(\int_{-\infty}^{\bar x}\exp\,(-2y^2)dy\biggr)^{-1}=-\infty.\cr}$$
Hence, for any $\ep\in (0,1)$ there exists ${\bar x}^*(\ep)<0$ such that
$$\exp\,(-2({\bar x}^*)^2)
\biggl(\int_{-\infty}^{{\bar x}^*}\exp\,(-2y^2)dy\biggr)^{-1}=
-{4{\bar x}^*(\ep)\over {1-\ep}}.\eqno(2.12)$$
Besides, ${\bar x}^*(\ep)\to 0$ when $\ep\to 1$. Let us take
$n\in [Ct+(2+2{\bar x}^*)\sqrt{Ct}-\Gamma,Ct+(2+2{\bar x}^*)\sqrt{Ct}+\Gamma]$.
Then ${{n-Ct-2\sqrt{Ct}}\over \sqrt{Ct}}=2{\bar x}^*+O(1/\sqrt{t})$.

We have now from one side
$$1-\v^{(-1)}\bigl({{n-2\sqrt{Ct}}\over t}\bigr)={C\over \v^{\prime}(1)}
{(-2{\bar x}^*)\over \sqrt{Ct}}+O(1/t).$$
From the other side we obtain using (2.12)
$$\eqalign{
&{1\over \sqrt{t}}\psi\bigl({{n-Ct-2\sqrt{Ct}}\over \sqrt{Ct}}\bigr)=
{1\over \sqrt{t}}\psi(2({\bar x}^*+O(1/\sqrt{t}))=\cr
&{C\over \v^{\prime}(1)\sqrt{t}}\exp\,(-2({\bar x}^*+O(1/\sqrt{t}))^2)
\biggl(\int_{-\infty}^{{\bar x}^*+O(1/\sqrt{t})}\exp\,(-2y^2)dy\biggr)^{-1}=
\cr
&{C\over \v^{\prime}(1)\sqrt{t}}\bigl({-4{\bar x}^*\over {1-\ep}}\bigr)+
O(1/t).\cr}$$
If ${2\over \sqrt{C}}<{4\over {1-\ep}}$ then there exists $t_0>0$ such that
$$1-\v^{(-1)}\bigl({{n-2\sqrt{Ct}}\over t}\bigr)<
{1\over \sqrt{t}}\psi\bigl({{n-Ct-2\sqrt{Ct}}\over \sqrt{Ct}}\bigr).$$
Besides, if $(1-\ep)$ small enough we have $-2{\bar x}^*\in [0,1]$.
So, we can finish the proof by putting $\delta=-2{\bar x}^*$.

\vskip 2 mm
{\bf Proof of Lemma 1.}
Let the function $\v(F)$ be extended for negative values of $F$ as a
smooth strictly decreasing function. Then there exists a wave-train solution
${\tilde F}_{\sigma}(n-C_{\sigma}t)$ for (2.1) with overfall
$(-\sigma,1)$, $\sigma>0$. Put $\sigma=\sigma(t)=\exp(-t^{1/3})$.
Proposition 1, Lemma 5, Lemma 6 from [7] together with Lemma 2, Lemma 3
above imply the following  statement.

For any $\delta\in (0,1)$, $l>1$, $A>2\sqrt{C}$ there exist $t_0>0$ and
increasing functions $\gamma_1(t)=O(t^{1/3})$, $\gamma_2(t)=2\sqrt{Clt}+a(l)$
such that
$$F^-(n,t)=\left\{\matrix{
{\tilde F}_{\sigma(t)}(n-Ct-\gamma_1),\ &n\le Ct+\sqrt{Clt}+a(l),
\hfill\cr
\v^{(-1)}\bigl({{n-\gamma_1-\gamma_2}\over t}\bigr),\ &
Ct+\sqrt{Clt}+a(l)<n<Ct+\gamma_1+\gamma_2-\delta\sqrt{Ct},\cr
1-{1\over \sqrt{t}}\psi\bigl({{n-Ct-\gamma_1-\gamma_2}\over \sqrt{Ct}}
\bigr),\ &Ct+\gamma_1+\gamma_2-\delta\sqrt{Ct}\le n<Ct+A\sqrt{Ct},\cr
1-\delta,\ &n\ge Ct+A\sqrt{t}\cr}\right.\eqno(2.13)$$
is a subsolution for (2.1), if $t\ge t_0$.

This statement and comparison principle from [6] imply that for any solution
$F(n,t)$ of the Cauchy problem (2.1), (2.2) there exists $T>0$ such that
$$F(n,t)>F^-(n,t+T),\eqno(2.14)$$
if $n\in\Z$, $t\ge -T+t_0$.

Lemma 1 follows from (2.13) and (2.14).

\vskip 2 mm
{\bf Proof of Proposition 1.}
Put $\kappa(t)=\{Ct+A\sqrt{t}\}$, $0\le\kappa(t)<1$, $N(t)=[Ct+A\sqrt{t}]$,
$F=F(N(t),t)$, $F_1=F(N(t)+1,t)$, ${\tilde F}={\tilde F}(N(t)-Ct+d_A(t))$,
${\tilde F}_1={\tilde F}(N(t)+1-Ct+d_A(t))$.
Proposition 3 from [7] implies the following asymptotic formula for
$$d_A^{\prime}(t)\buildrel \rm def \over = {d\over dt}(d_A(t))$$
$$\eqalign{
&d_A^{\prime}(t)(1+O(1/\sqrt{t}))=C(1-\kappa)(F_1-F)-C(1-\kappa)
({\tilde F}_1-{\tilde F})+\cr
&{A(1-{\tilde F}_1)\over 2\sqrt{t}}+{1\over 2}\v^{\prime}(1)
(1-{\tilde F}_1)^2-\bigl(
{A(1-F_1)\over 2\sqrt{t}}+{1\over 2}\v^{\prime}(1)
(1-F_1)^2\bigr).\cr}\eqno(2.15)$$
Let us estimate now all terms of (2.15). The assumption
$\Delta F(n,t_0)\ge 0$ implies (by Theorem 1 in [5]) that
$\Delta F(n,t_0)\ge 0\ \ \forall\ t\ge t_0$. From this and from inequality
(2.7) it follows for $n\ge N(t)+1$ and $t\ge t_0$:
$$\eqalign{
&0\le 1-F(n,t)\le 1-F_1\le 1-F(N(t)+1,t+T)\le {1\over \sqrt{t}}
\psi(A-\delta_0-2)\le\cr
&{C\over \v^{\prime}(1)\sqrt{t}}\exp\,(-(A-\delta_0-2)^2/2)
\biggl(\int_{-\infty}^0\exp\,(-2y^2)dy\biggr)^{-1}\le\cr
&O\bigl({1\over \sqrt{t}}\exp\,(-(A-\delta_0-2)^2/2)\bigr).\cr}\eqno(2.16)$$
From (2.16) and from inequality (1.12) of Theorem 2ii) for $t\ge t_0\ge A^2$,
$n\ge N(t)$
we obtain the crucial inequality
$$F_1-F=O\bigl({A\over t}\exp\,(-(A-\delta_0-2)^2/2)\bigr).\eqno(2.17)$$
From [5] (Theorems 2, $2^{\prime}$) and [6] (Theorems 6.1, 6.2) it follows
asymptotic formula
$${\tilde F}_1=1-{C\over \v^{\prime}(1)(A\sqrt{t}+d_A(t))}+
O\bigl({1\over (A\sqrt{t}+d_A(t))^2}\bigr).$$
This formula and estimate $d_A(t)\ge 0$ (see (2.4)) gives inequalities
$$\eqalign{
&0<1-{\tilde F}_1\le O\bigl({1\over A\sqrt{t}}\bigr),\cr
&{\tilde F}_1-{\tilde F}=O\bigl({1\over A^2t}\bigr).\cr}\eqno(2.18)$$

Let us put estimates (2.16)-(2.18) into formula (2.15). We obtain
$$\eqalign{
&d_A^{\prime}(t)(1+O(1/\sqrt{t}))=C(1-\kappa)(F_1-F)-C(1-\kappa)
({\tilde F}_1-{\tilde F})+\cr
&{C\over 2\v^{\prime}(1)t}+{1\over 2}{C^2\over \v^{\prime}(1)A^2t}-
{(1-F_1)\over 2}\bigl({A\over \sqrt{t}}+\v^{\prime}(1)(1-F_1)\bigr)=\cr
&(1-\kappa)O\bigl({A\over t}\exp\,(-(A-\delta_0-2)^2/2)\bigr)-
(1-\kappa)O\bigl({1\over A^2t}\bigr)+\cr
&{C\over 2\v^{\prime}(1)t}+{1\over 2}{C^2\over \v^{\prime}(1)A^2t}+
O\bigl({A\over t}\exp\,(-(A-\delta_0-2)^2/2)\bigr)=\cr
&{C\over 2\v^{\prime}(1)t}+O\bigl({1\over A^2t}\bigr)+
O\bigl({A\over t}\exp\,(-(A-\delta_0-2)^2/2)\bigr).\cr}\eqno(2.19)$$
Estimate (2.9) implies asymptotic formula
$$d_A(t)={C\over 2\v^{\prime}(1)}\ln\,t+O(1/A^2)\ln\,t + const.$$
From result (2.5) it follows that for any $A_1>2\sqrt{C}$ and $A_2>2\sqrt{C}$
we have $d_{A_1}(t)-d_{A_2}(t)\to 0$, $t\to\infty$.

Hence,
$$d_A(t)={C\over 2\v^{\prime}(1)}\ln\,t + const +o(1).$$

\vskip 4 mm
{\bf 3. A priori estimates for local solutions of Burgers type equations.}
\vskip 2 mm
Without loss of generality we will put further $C=1$ and $\ep=1$.
Otherwise we make substitutions: $t\to Ct/\ep$, $x\to x/\ep$ for the
equation (1.2) and $t\to C^2t/\ep$, $x\to Cx/\ep$ for the equation (1.1). We
will give here a complete proof of Theorem 2ii) which is sufficient
for all current applications and a sketch of the proof of
Theorem $2ii)^{\prime}$. Theorem 2i) will be proved in a separate paper.

The first step in the proof of Theorem 2ii) is the Green-Poisson type
representation formula (for function $u$ in $\Omega_{\sigma}$) associated
with operator $u\mapsto u_t^{\prime}+\Delta u$, where
$\Delta u\buildrel \rm def \over = u(x,t)-u(x-1,t)$,
$u_t^{\prime}={\pa u(x,t)\over \pa t}$.

Let $\chi_0\ \ :\ \ \R\to\R$ be a smooth cut-off function  such that
$$\eqalign{
&0\le\chi_0\le 1,\ \chi_0\bigg|_{(-\infty,a_1)}\equiv 0,\
\chi_0\bigg|_{[{\tilde a}_1,+\infty]}\equiv 1,\ 0<a_1<{\tilde a}_1<\infty,\cr
&|\chi_{0t}^{\prime}|\le {A_0\over \delta}\ \ {\rm and}\ \
|\chi_{0t}^{\prime\prime}|\le {A_0\over \delta^2},\cr}\eqno(3.1)$$
 where $\delta={\tilde a}_1-a_1$.
Put $\chi(x,t)=\chi_0\bigl({{x-t}\over \sqrt{t}}\bigr)$.

\vskip 2 mm
{\bf Proposition 2.}
{\it Let function} $u(x,t)$ {\it be defined in the domain}
$\Omega_{\sigma}=\{(x,t)\ \ :\ \ a_1<\bar x\buildrel \rm def \over =
{{x-t}\over \sqrt{t}}<a_2+\sigma\sqrt{t}\}$,
$\sigma>0$ {\it and} ${\tilde u}(x,t)=u(x,t)\cdot \chi(x,t)$. {\it Let}
$0<\sigma_0<\sigma$ {\it and} $\alpha\in \bigl({{1+\sigma_0}\over {1+\sigma}},
1\bigr)$. {\it Then function} ${\tilde u}$ {\it can be represented in}
$\Omega_{\sigma_0}$ {\it by the following formula of the Green-Poisson type}
$$\eqalign{
&{\tilde u}(x,t)=\int_{-\infty}^{\infty}G(x-\xi,t-\alpha t){\tilde u}
(\xi,\alpha t)d\xi+\cr
&\int_{\alpha t}^td\tau\int_{-\infty}^{\infty}G(x-\xi,t-\tau)
({\tilde u}_{\tau}^{\prime}+\Delta {\tilde u})(\xi,\tau)d\xi,\cr}\eqno(3.2)$$
{\it where}
$$G(x,t)={1\over 2\pi}\int_{-\infty}^{\infty}\exp\,(-i\xi x)
\exp\,([e^{i\xi}-1]t)\,d\xi.$$
{\it Besides},
$$G(x,t)=\sum_{n=-\infty}^{\infty}G_n(t)\delta(n-x),\eqno(3.3)$$
{\it where}
$$\left\{\matrix{
G_n(t)=0,\ \ &{\it if}\ \ n<0,\hfill\cr
G_n(t)={t^n\over n!}e^{-t},\ \ &{\it if}\ \ n\ge 0 \cr}\right.$$
is Poisson-distribution.

This statement is certainly classical but we did not find the precise
reference. So, we will indicate the abridge proof.

The operator ${\pa\over \pa t}+\Delta$ can be considered as a parabolic
operator of infinite order in $x$ and it can be represented by the
following formula
$${\pa\over \pa t}+\Delta={\pa\over \pa t}+(1-\exp\,(-{\pa\over \pa x})).$$
We will apply further to the Cauchy problem for this operator the same
Fourier method as for parabolic operator of finite order and we will obtain
(3.2). The formula (3.3) is the Fourier inversion formula for the
classical Poisson distribution through its caracteristic function.

It is important to remark that the function ${\tilde u}(\xi,\tau)$ is well defined
for $(\xi,\tau)\ \ :\ \ \xi<\tau +a_2\sqrt{\tau}+\sigma\tau$,
${\tilde u}(\xi,\tau)\equiv 0$ for $\xi\le\tau+a_1\sqrt{\tau}$ and
function $\xi\mapsto G(x-\xi,t-\tau)$ is equal to zero for
$\xi>x=t+{\bar x}\sqrt{t}$.
So, the function $\xi\mapsto {\tilde u}(\xi,\tau)\cdot G(x-\xi,t-\tau)$ can
be naturally interpreted in the formula (3.2) as a function with
compact support in $\R$ if the following inequality is satisfied
$$\eqalign{
&\tau+a_2\sqrt{\tau}+\sigma\tau\ge
x=t+{\bar x}\sqrt{t}\ \ {\rm for}\cr
&{\bar x}\in (a_1,a_2+\sigma_0\sqrt{t}),\ \ \sigma_0<\sigma\ \ {\rm and}\ \
\tau\ge\alpha t\ge t_0(\sigma,\sigma_0).\cr}$$
In order to satisfy these inequalities we choose $\alpha\in (0,1)$ such
that for $t>t_0(\sigma,\sigma_0)$ the following inequality is valid
$$\alpha t+a_2\sqrt{\alpha t}+\sigma\alpha t>t+a_2\sqrt{t}+\sigma_0t,$$
i.e. we must take $\alpha>{{1+\sigma_0}\over {1+\sigma}}$.

\vskip 2 mm
{\bf Corollary}\ (Integral representation for $\Delta u(x,t)$).
{\it Let function} $u(x,t)$ {\it satisfy} (1.2) {\it in} $\Omega_{\sigma}$
{\it with} $\v(0)=C=1$ {\it and} $\ep=1$. {\it Put} $\v_0=\v-C$.  {\it Then
in assumption of Proposition 2 for}
$$(x,t)\in{\tilde\Omega}_{\sigma_0}=\{(x,t)\in\Omega_{\sigma_0}\ :\ x\ge t+
{\tilde a}_1\sqrt{t}\},\ \ \sigma_0<\sigma,\ \ t>t^*=\alpha t\ge t_0,\ \
\alpha\in ({{1+\sigma_0}\over {1+\sigma}},1)$$
{\it we have the equality}
$$\Delta u(x,t)=I_0u+I_1u+I_2u+I_3u+I_4u,\eqno(3.4)$$
{\it where}
$$\eqalign{
&I_0u(x,t)=-\int_{t^*}^td\tau\int_{{\bar \xi}>{\tilde a}_1}
\Delta_xG(x-\xi,t-\tau)\v_0(u)\Delta u(\xi,\tau)d\xi,\cr
&I_1u(x,t)=-\int_{t^*}^td\tau\int_{{\bar \xi}\in [a_1,{\tilde a}_1]}
\Delta_xG(x-\xi,t-\tau)\v_0(u)\Delta u(\xi,\tau)\chi(\xi,\tau)d\xi,\cr
&I_2u(x,t)=\int_{{\bar\xi}\ge a_1}\Delta_xG(x-\xi,t-t^*)u(\xi,t^*)
\chi(\xi,t^*)d\xi,\cr
&I_3u(x,t)=\int_{t^*}^td\tau\int_{{\bar \xi}\in [a_1,{\tilde a}_1]}
\Delta_xG(x-\xi,t-\tau)(u\chi_{\tau}^{\prime}+u\Delta\chi)(\xi,\tau)d\xi,\cr
&I_4u(x,t)=-\int_{t^*}^td\tau\int_{{\bar \xi}\in [a_1,{\tilde a}_1]}
\Delta_xG(x-\xi,t-\tau)\Delta u(\xi,\tau)\Delta\chi(\xi,\tau)d\xi.\cr}$$

\vskip 2 mm
{\bf Remark.}
We will use below several times the following simple relation:\ let
$u=u(x)$, $v=v(x)$, then $\Delta (u\cdot v)=u\cdot\Delta v + v(x-1)\Delta u$,
where $\Delta u\buildrel \rm def \over = u(x)-u(x-1)$.

\vskip 2 mm
{\bf Proof of Corollary.}
We have relations
$$\eqalign{
&{\tilde u}(\xi,\tau)=u(\xi,\tau)\cdot\chi(\xi,\tau),\cr
&{\tilde u}_{\tau}^{\prime}=(u\cdot\chi)_{\tau}^{\prime}=u_{\tau}^{\prime}\cdot
\chi + u\chi_{\tau}^{\prime},\cr
&\Delta {\tilde u}=\Delta(u\cdot\chi)=\Delta u\cdot\chi(\xi-1,t)+
u\cdot\Delta\chi=\Delta u\cdot\chi+u(\xi-1,t)\Delta\chi.\cr}$$
Using (1.2) we obtain
$$(u_{\tau}^{\prime}+\Delta u)\cdot\chi=-\v_0(u)\Delta u\chi=-\v_0(u)
(\Delta {\tilde u}-u(\xi-1,t)\cdot\Delta\chi),$$
$$\eqalign{
&({\tilde u}_{\tau}^{\prime}+\Delta {\tilde u})=-\v_0(u)\Delta {\tilde u}+
\v_0(u)\cdot u(\xi-1,\tau)\Delta\chi+u(\chi_{\tau}^{\prime}+\Delta\chi)-
\Delta u\cdot\Delta\chi=\cr
&-\v_0(u)\Delta u\cdot\chi + u\cdot(\chi_{\tau}^{\prime}+\Delta\chi)
-\Delta u\cdot\Delta\chi.\cr}$$

Plugging these relations into (3.2) and using the equality
${\tilde u}(\xi,\tau)=u(\xi,\tau)$ for $\bar\xi>{\tilde a}_1$ we obtain (3.4).

For the estimates of terms $I_1u$, $I_2u$, $I_3u$, $I_4u$ in  formula (3.4)
we will use elementary estimates for cut-off function $\chi(x,t)$ and
rather precise estimates for Green-Poisson function $G(x,t)$.

\vskip 2 mm
{\bf Lemma 4.}
{\it Let} $\chi(x,t)$ {\it be cut-off function defined by} (3.1). {\it Then
the following estimates for derivatives of} $\chi$ {\it are valid}
$$\eqalign{
&|\Delta\chi(x,t)|\le {A_0\over \delta\sqrt{t}},\ \
|\Delta^2\chi(x,t)|\le {A_0\over \delta^2t}\ \ {\it and}\cr
&|(\chi_t^{\prime}+\Delta\chi)(x,t)|\le {A_0\over t}\bigl({1\over \delta^2}+
{{\tilde a}_1\over 2\delta}\bigr),\cr}$$
{\it where} $(x,t)\in\Omega_{\sigma}$, $\delta={\tilde a}_1-a_1$.

\vskip 2 mm
{\bf Proof.}
We have

$$\eqalign{
&\chi^{\prime}(x,t)=-\bigl({1\over \sqrt{t}}+{{x-t}\over 2t^{3/2}}\bigr)
\chi_0^{\prime}\bigl({{x-t}\over \sqrt{t}}\bigr);\cr
&\Delta\chi(x,t)={1\over \sqrt{t}}\int_{x-1}^x\chi_0^{\prime}
\bigl({{y-t}\over \sqrt{t}}\bigr)dy;\cr}$$
$$\eqalign{
&(\chi^{\prime}+\Delta\chi)(x,t)=-{1\over \sqrt{t}}
\int_{x-1}^x\bigl(\chi_0^{\prime}\bigl({{x-t}\over \sqrt{t}}\bigr)-
\chi_0^{\prime}\bigl({{y-t}\over \sqrt{t}}\bigr)\bigr)dy-
\chi_0^{\prime}\bigl({{x-t}\over \sqrt{t}}\bigr){{x-t}\over 2t^{3/2}}=\cr
&-{1\over t}\int_{x-1}^x\int_y^x\chi_0^{\prime\prime}
\bigl({{z-t}\over \sqrt{t}}\bigr)dzdy-
\chi_0^{\prime}\bigl({{x-t}\over \sqrt{t}}\bigr){{x-t}\over 2t^{3/2}}.\cr}$$
From these relations and from estimates (3.1) for $\chi_0$ we obtain
necessary estimates for $\chi(x,t)$.

\vskip 2 mm
{\bf Lemma 5.} (Estimates for Green-Poisson distribution $G(x,t)$).
{\it Let}

\noindent
$G(x,t)=\sum_{n=0}^{\infty}G_n(t)\delta(n-x)$ {\it be the
Poisson distribution} (3.3). {\it The following estimates for}
$\{G_n(t)\}$ {\it are valid}

i) {\it if} $p=n-t\ge 0$ {\it then}
$$G_n(t)\le {1\over \sqrt{2\pi n}}e^{-p^2/(2n)}.$$

ii) {\it if} $q=t-n>0$, $q\le t$ {\it then}
$$G_n(t)\le {1\over \sqrt{2\pi n}}e^{-q^2/(2t)}.$$

iii) {\it if} $n=t+a\sqrt{t}$ {\it then}
$$\eqalign{
&G_n(t)={1\over \sqrt{2\pi n}}\exp\,\bigl(-{(n-t)^2\over 2t}\bigr)
\biggl(1+O\biggl({(n-t)^3\over t^2}\biggr)\biggr)=\cr
&{1\over \sqrt{2\pi t}}\exp\,\bigl(-{a^2\over 2}\bigr)
\biggl(1+O\biggl({{a+a^3}\over \sqrt{t}}\biggr)\biggr).\cr}$$

\vskip 2 mm
{\bf Proof.}
By Stirling's formula we have
$$n!=\sqrt{2\pi n}\bigl({n\over e}\bigr)^n\bigl(1+O\bigl({1\over n}\bigr)
\bigr).$$
Then
$$G_n(t)={1\over \sqrt{2\pi n}}e^{n\ln\,t-n\ln\,n+n-t}
\bigl(1-O\bigl({1\over n}\bigr)\bigr).$$
If $p=n-t>0$ then
$$\ln\,{t\over n}=\ln\,\bigl(1-{p\over n}\bigr)=
-{p\over n}-{p^2\over 2n^2}-\ldots.$$
If $q=t-n>0$ then
$$\ln\,{n\over t}=\ln\,\bigl(1-{q\over t}\bigr)=
-{q\over t}-{q^2\over 2t^2}-\ldots.$$
Hence,
$$G_n(t)={1\over \sqrt{2\pi n}}e^{-p^2/(2n)}
\bigl(1-O\bigl({1\over n}\bigr)\bigr),\ \ {\rm if}\ \ p=n-t>0$$
and
$$\eqalign{
&G_n(t)={1\over \sqrt{2\pi n}}e^{-q^2/(2t)-(1/2-1/3)(q^3/t^2)-\ldots}
\bigl(1-O\bigl({1\over n}\bigr)\bigr),\cr
&{\rm if}\ \ q=t-n>0,\ \ q<t.\cr}$$
These relations give i), ii) and iii).

\vskip 2 mm
{\bf Lemma 6.} (Estimates for $\Delta G(x,t)$).
{\it Let} $G(x,t)=\sum_{n=0}^{\infty}G_n(t)\delta(n-x)$ {\it be the
Poisson distribution. We put}
$$\eqalign{
&\Delta G_n(t)=G_n(t)-G_{n-1}(t),\cr
&\Delta_xG(x-\xi,t-\tau)=G(x-\xi,t-\tau)-G(x-1-\xi,t-\tau),\cr
&\bar\xi={{\xi-\tau}\over \sqrt{\tau}},\ \ \bar x={{x-t}\over \sqrt{t}}.\cr}
$$
{\it Then the following estimates are valid}

i)
$$\eqalign{
&\Delta G_n(t)=G_n(t){(t-n)\over t},\ \ {\it and\ as\ consequence}\cr
&\Delta G_n(t)>0,\ \ {\it if}\ \ n<t,\ \ \Delta G_n(t)<0,\ \ {\it if}\ \
n>t;\cr
&\Delta^2G_n(t)=G_n(t)\bigl(1-{2n\over t}+{n(n-1)\over t^2}\bigr),\ \
{\it and\ as\ consequence}\cr
&\Delta^2G_n(t)<0,\ \ {\it if}\ \ n-t-{1\over 2}\in
\bigl(-\sqrt{t+1/4},+\sqrt{t+1/4}\bigr),\cr
&\Delta^2G_n(t)\ge 0,\ \ {\it if}\ \ n-t-{1\over 2}\notin
\bigl(-\sqrt{t+1/4},+\sqrt{t+1/4}\bigr);\cr}$$

ii) $\forall\ s\ge 0$ {\it and} $p\ge 0$ {\it we have inequalities}
$$\eqalign{
&-\Delta G_{p+s}(s)\le A_1s^{-3/2}p\,\exp\,\bigl(-{p^2\over 4s}\bigr),\ \
{\it if}\ \ p<s,\cr
&-\Delta G_{p+s}(s)\le A_1p^{-1/2}e^{-p/4},\ \ {\it if}\ \ p>s;\cr}$$

iii) $\forall\ s\ge 0$ {\it and} $q\in (0,s)$ {\it we have inequalities}
$$\Delta G_{s-q}(s)\le A_1{q\over s\sqrt{s-q}}\exp\,\bigl(-{q^2\over 2s}
\bigr);$$

iv)
$$\eqalign{
&\sum_{n=-1}^{\infty}|\Delta G_n(t)|=\min\,\bigl\{2,{2\over \sqrt{2\pi t}}
\bigl(1+O\bigl({1\over \sqrt{t}}\bigr)\bigr)\bigr\};\cr
&\sum_{n=-2}^{\infty}|\Delta^2G_n(t)|=\min\,\bigl\{4,{4\over \sqrt{2\pi e}}
{1\over t}\bigl(1+O\bigl({1\over \sqrt{t}}\bigr)\bigr)\bigr\};\cr}$$

v) $\forall\ \bar x>{\tilde a}_1$ {\it and} $t>\tau>\alpha t$ {\it we
have inequality}
$$\eqalign{
&I=\int_{\bar\xi>{\tilde a}_1}|\Delta_xG(x-\xi,t-\tau)|
\bigl(1+\ln_+{1\over {\bar\xi-{\tilde a}_1}}\bigr)(1+\bar\xi)d\xi\le\cr
&{A_1\over \sqrt{t-\tau}} \bigl(1+\sqrt{(1-\alpha)/\alpha}\bigr)
\bigl(1+\ln_+{1\over {\bar x-{\tilde a}_1}}\bigr)(1+{\bar x}/\sqrt{\alpha}),
\cr}$$
{\it where} $A_1$ {\it is absolute constant}.

\vskip 2 mm
{\bf Remark.}
We will use further several times the differential and integral relations:
$$\eqalign{
&-\Delta_{\xi}(G(x-\xi-1,t-\tau)u(\xi))=G(x-\xi,t-\tau)\Delta u(\xi)+
\Delta_{\xi}G(x-\xi-1,t-\tau)\cdot u(\xi);\cr
&-\Delta_{\xi}G(x-\xi-1,t-\tau)=\Delta_xG(x-\xi,t-\tau);\cr}$$
if $G(x-\xi-1,t-\tau)\cdot u(\xi)$ has compact support with respect to $\xi$
then
$$\eqalign{
&-\int_{\xi\in\R}\Delta_{\xi}(G(x-\xi-1,t-\tau)\cdot u(\xi))d\xi=0\ \
{\rm and\ hence}\cr
&-\int_{\xi\in\R}\Delta_xG(x-\xi,t-\tau)\cdot u(\xi)d\xi=
\int_{\xi\in\R}G(x-\xi,t-\tau)\Delta u(\xi)d\xi.\cr}$$

\vskip 2 mm
{\bf Proof of Lemma 6.}

i) We have from (3.3)
$$\eqalign{
&\Delta G_n(t)=\bigl({t^n\over n!}-{t^{n-1}\over (n-1)!}\bigr)e^{-t}=
G_n(t){(t-n)\over t},\cr
&\Delta^2G_n(t)=\bigl({t^n\over n!}-2{t^{n-1}\over (n-1)!}+{t^{n-2}\over
(n-2)!}\bigr)e^{-t}=G_n(t)\bigl(1-{2n\over t}+{n(n-1)\over t^2}\bigr);\cr}$$

ii) follows from i) and Lemma 5i).

iii) follows from i) and Lemma 5ii).

iv) Putting in i) $p=n-t=a\sqrt{t}$ and using Lemma 5iii) we obtain
$$\eqalign{
&\Delta G_n(t)={1\over \sqrt{2\pi t}}e^{-a^2/2}\bigl(-{a\over \sqrt{t}}\bigr)
\bigl(1-O\bigl({a^3\over \sqrt{t}}\bigr)\bigr)\ \ {\rm and}\cr
&{\rm as\ consequence}\ \  \Delta G_n(t)=0\ \ {\rm if}\ \ a=0.\cr}$$

So,
$$\sum_{n=-1}^{\infty}|\Delta G_n(t)|=\bigl(\sum_{n\ge t}\Delta G_n(t)-
\sum_{n\le t}\Delta G_n(t)\bigr).$$
Then
$$\eqalign{
&\sum_{n=-1}^{\infty}|\Delta G_n(t)|=[G_{[t]}(t)-G_{-\infty}(t)]-
[G_{+\infty}(t)-G_{[t]}(t)]=\cr
&2G_{[t]}(t)={2\over \sqrt{2\pi t}}\bigl(1+O\bigl({1\over \sqrt{t}}\bigr)
\bigr),\ \ t\ge t_0.\cr}$$
For all $t>t_0$ we have
$$\sum_{n=-1}^{\infty}|\Delta G_n(t)|=\min\,\bigl\{2,{2\over \sqrt{2\pi t}}
\bigl(1+O\bigl({1\over \sqrt{t}}\bigr)\bigr)\bigr\}.$$
By similar arguments we have
$$\eqalign{
&\sum_{n=-2}^{\infty}|\Delta^2G_n(t)|=\sum_{-\infty}^{[t-\sqrt{t}]}
\Delta^2G_n(t)-\sum_{[t-\sqrt{t}]+1}^{[t+\sqrt{t}]}\Delta^2G_n(t)+
\sum_{[t+\sqrt{t}]+1}^{\infty}\Delta^2G_n(t)=\cr
&2\Delta G_{[t-\sqrt{t}]}(t)+2|\Delta G_{[t+\sqrt{t}]}(t)|=
{4\over t\sqrt{2\pi e}}\bigl(1+O\bigl({1\over \sqrt{t}}\bigr)\bigr),\ \
t\ge t_0.\cr}$$
For all $t>0$ we have
$$\sum_{n=-2}^{\infty}|\Delta^2G_n(t)|=\min\,\bigl\{4,{4\over t\sqrt{2\pi e}}
\bigl(1+O\bigl({1\over \sqrt{t}}\bigr)\bigr)\bigr\}.$$

v) Put $t-\tau=s$, $x-\xi=y$. We have $p=y-s=\bar x\sqrt{t}-
\bar\xi\sqrt{\tau}$. Put $I=I_++I_-$, where
$$I_{\pm}=\int_{\pm\Delta_xG<0}|\Delta_xG(x-\xi,t-\tau)|
\bigl(1+\ln_+{1\over (\bar\xi-{\tilde a}_1)}\bigr)(1+\bar\xi)
d\xi.$$
By part i) $\Delta_xG(x-\xi,t-\tau)<0$ iff $p=(x-\xi)-(t-\tau)>0$.

Hence, $I_+=I_+^{\prime}+I_+^{\prime\prime}$, where
$$\eqalign{
&I_+^{\prime}=-\int_{\bar\xi>{\tilde a}_1\ :\ 0<p<s}\Delta_xG(x-\xi,t-\tau)
\bigl(1+\ln_+{1\over (\bar\xi-{\tilde a}_1)}\bigr)(1+\bar\xi)
d\xi,\cr
&I_+^{\prime\prime}=-\int_{\bar\xi>{\tilde a}_1\ :\ p>s}
\Delta_xG(x-\xi,t-\tau)
\bigl(1+\ln_+{1\over (\bar\xi-{\tilde a}_1)}\bigr)(1+\bar\xi)
d\xi.\cr}$$
Put $p_1=\bar x\sqrt{t}-{\tilde a}_1\sqrt{\tau}$. We have
${\tilde a}_1-\bar\xi={{p-p_1}\over \sqrt{\tau}}<0$ and
$p=p_1$ iff\ \ $\bar\xi={\tilde a}_1$.

For $I_+^{\prime}$ when $p\in (0,s)$ we use ii) and obtain
$$\eqalign{
&I_+^{\prime}\le As^{-3/2}\int_0^{p_1}e^{-p^2/(4s)}p\,\bigl(1+
\ln_+{\sqrt{\tau}\over {p_1-p}}\bigr)\bigl(1+{\tilde a}_1+{{p_1-p}\over
\sqrt{\tau}}\bigr)dp\le  \cr
&({\rm putting}\ \ p=\rho\sqrt{s})\cr
&As^{-1/2}\int_0^{p_1/\sqrt{s}}e^{-\rho^2/4}\rho\,\bigl(1+
\ln_+{\sqrt{\tau/s}\over (p_1/\sqrt{s}-\rho)}\bigr)
\bigl(1+{\tilde a}_1+{p_1\over \sqrt{\tau}}-\rho\bigr)
d\rho\le \cr
&({\rm by\ Lemma}\ A_1\ \ {\rm of\ Appendix})\cr
&As^{-1/2}\bigl(1+\ln_+{\sqrt{\tau}\over p_1}\bigr)
\bigl(1+{\tilde a}_1+{p_1\over \sqrt{\tau}}\bigr)\le
As^{-1/2}\bigl(1+\ln_+{1\over {\bar x\sqrt{t/\tau}-{\tilde a}_1}}\bigr)
(1+{\bar x}\sqrt{t/\tau})\le\cr
&A_1s^{-1/2}\bigl(1+\ln_+{1\over (\bar x-{\tilde a}_1)}\bigr)
(1+{\bar x}/\sqrt{\alpha}).\cr}$$
For $I_+^{\prime\prime}$ when $p>s$ we use ii) and obtain
$$\eqalign{
&I_+^{\prime\prime}\le As^{-1/2}\int_0^{p_1}e^{-p/4}\bigl(1+
\ln_+{\sqrt{\tau}\over {p_1-p}}\bigr)
\bigl(1+{\tilde a}_1+{p_1\over \sqrt{\tau}}-\rho\bigr)dp
\le\cr
&As^{-1/2}\bigl(1+\ln_+{\sqrt{\tau}\over p_1}\bigr)
(1+{\bar x}\sqrt{t/\tau})\le
A_1s^{-1/2}\bigl(1+\ln_+{1\over (\bar x-{\tilde a}_1)}\bigr)
(1+{\bar x}/\sqrt{\alpha}).\cr}$$
Let us estimate now integral $I_-$. Put $q=(t-\tau)-(x-\xi)$. By part i)
$\Delta_xG(x-\xi,t-\tau)>0$ iff $q\in (0,s)$. We use now part iii) and
obtain
$$\eqalign{
&I_-\le {A\over s}\int_0^se^{-q^2/(2s)}
{q\over \sqrt{s-q}}
\bigl(1+\ln_+{\sqrt{\tau}\over {p_1+q}}\bigr)
\bigl(1+{\tilde a}_1+{{p_1+q}\over \sqrt{\tau}}\bigr)dq\le\cr
&{A\over s}\bigl(1+\ln_+{1\over {\bar x-{\tilde a}_1}}\bigr)
\biggl((1+{\bar x}/\sqrt{\alpha})\int_0^s{e^{-q^2/(2s)}q\over \sqrt{s-q}}dq+
{1\over \sqrt{\tau}}\int_0^s{e^{-q^2/(2s)}q^2\over \sqrt{s-q}}dq
\biggr),\cr}$$
where
$$\eqalign{
&\int_0^s{e^{-q^2/(2s)}q\over \sqrt{s-q}}dq\le
\int_0^{s/2}{e^{-q^2/(2s)}q\over \sqrt{s/2}}dq+
\int_{s/2}^s{e^{-s/8}s\over \sqrt{s-q}}dq\le\cr
&\sqrt{2s}\int_0^{s/8}e^{-y}dy+2se^{-s/8}\sqrt{s/2}=\sqrt{2s}
\bigl(1-e^{-s/2}+se^{-s/8}\bigr)=O(\sqrt{s}),\cr
&\int_0^s{e^{-q^2/(2s)}q^2\over \sqrt{s-q}}dq\le
\int_0^{s/2}{e^{-q^2/(2s)}q^2\over \sqrt{s/2}}dq+
\int_{s/2}^s{e^{-s/8}s^2\over \sqrt{s-q}}dq\le\cr
&2s\cdot\int_0^{s/8}\sqrt{y}e^{-y}dy+2s^2e^{-s/8}\sqrt{s/2}=
O(s).\cr}$$

Hence,
$$I_-\le {A_2\over \sqrt{s}}\bigl(1+\ln_+{1\over {\bar x-{\tilde a}_1}}
\bigr)
\bigl(1+\sqrt{(1-\alpha)/\alpha}\bigr)(1+{\bar x}/\sqrt{\alpha}).$$
Lemma 6 is proved.

Now we are ready to estimate terms $I_2u$ and $I_3u$ of formula (3.4).

\vskip 2 mm
{\bf Lemma 7.}
{\it Let function} $F=u$ {\it satisfy the conditions of Theorem} 2ii)
{\it and}
$\Delta u$ {\it is represented in} $\Omega_{\sigma}$ {\it by formula} (3.4),
$\alpha\ge\sup\{1/2, {{1+\sigma_0}\over {1+\sigma}}\}$, $\sigma_0<\sigma$.
{\it Then
terms} $I_2u$ {\it and} $I_3u$ {\it of formula} (3.4) {\it admit the
following estimates}
$$
|I_2u(x,t)|\le
A_2{\Gamma\cdot\bar x\over \sqrt{(1-\alpha)}t},\eqno(3.5)$$
$$\eqalign{
&|I_3u(x,t)|\le A_0{\Gamma\cdot {\tilde a}_1\over t^{3/2}}
\bigl({1\over \delta^2}+{{\tilde a}_1\over 2\delta}\bigr)
K^+,\cr
&{\it where}\ \ K^+=\int_{\alpha t}^td\tau\int_{\bar\xi\in [a_1,{\tilde a}_1]}
|\Delta_xG(x-\xi,t-\tau)|d\xi,\cr}\eqno(3.6)$$
$A_2$ {\it is absolute constant}, $\bar x\in ({\tilde a}_1,a_2+
\sigma_0\sqrt{t})$, $t>t_0(\sigma_0,\sigma)$.

\vskip 2 mm
{\bf Remark.}
$I_2u$ is the only term in representation (3.4), where $(1-\alpha)$ is in
the denominator.

\vskip 2 mm
{\bf Proof.}
The definitions of $I_2u$ and $I_3u$, condition (1.11) and Lemma 4 imply
estimates:
$$|I_2u(x,t)|\le\Gamma\int_{\bar\xi>a_1}|\Delta_xG(x-\xi,t-\alpha t)|
{\bar\xi d\xi\over \sqrt{\alpha t}},\eqno(3.7)$$
where $\bar\xi={{\xi-\alpha t}\over \sqrt{\alpha t}}$,
$$|I_3u(x,t)|\le\Gamma A_0\bigl({1\over \delta^2}+{{\tilde a}_1\over
2\delta}\bigr)
\int_{\alpha t}^td\tau\int_{\bar\xi\in (a_1,{\tilde a}_1)}
{|\Delta_xG(x-\xi,t-\tau)|\over \tau}
{\bar\xi d\xi\over \sqrt{\tau}},\eqno(3.8)$$
where $\bar\xi={{\xi-\tau}\over \sqrt{\tau}}$.

Using Lemmas 6i), 5i), 6iv) (see also (3.14))
we obtain further from (3.7)
$$\eqalign{
&|I_2u(x,t)|\le\cr
&{\Gamma\over \sqrt{\alpha t}}
\biggl[\int_{\bar\xi<\bar x\sqrt{t/(\alpha t)}}
\Delta_{\xi}G(x-\xi,t-\alpha t)\bar\xi d\xi-
\int_{\bar\xi>\bar x\sqrt{t/(\alpha t)}}
\Delta_{\xi}G(x-\xi,t-\alpha t)\bar\xi d\xi\biggr]\le\cr
&{\Gamma\over \sqrt{\alpha t}}
\biggl[-\int_{\bar\xi<\bar x/\sqrt{\alpha}}G\cdot\Delta_{\xi}\bar\xi d\xi+
\int_{\bar\xi>\bar x/\sqrt{\alpha}}G\cdot\Delta_{\xi}\bar\xi d\xi+
G\bar\xi\big|_{a_1}^{\bar x/\sqrt{\alpha}}-
G\bar\xi\big|_{\bar x/\sqrt{\alpha}}^{\bar x/\sqrt{\alpha}+
(1-\alpha)\sqrt{t}/\sqrt{\alpha}}\biggr]\le \cr
&{\Gamma\over \sqrt{\alpha t}}\biggl[\int_{\bar\xi>\bar x/\sqrt{\alpha}}
{1\over \sqrt{\alpha t}}G(x-\xi,t-\alpha t)d\xi +
2G(t-\alpha t,t-\alpha t){\bar x\over \sqrt{\alpha}}\biggr]\le\cr
&{\Gamma\over \sqrt{\alpha t}}
\biggl({1\over \sqrt{\alpha t}}
+
{2\over \sqrt{2\pi(t-\alpha t)}}\bigl({\bar x\over \sqrt{\alpha}}\bigr)
\biggr)\le\cr
&{A_2\Gamma\over \sqrt{(1-\alpha)\alpha}}{1\over t}\bigl({\bar x\over
\sqrt{\alpha}}\bigr),\ \ {\rm if}\ \ t\ge t_0.\cr}$$

From (3.8) we deduce
$$|I_3u(x,t)|\le {\Gamma\,A_0\over t^{3/2}}\biggl({1\over \delta^2}+
{{\tilde a}_1\over 2\delta}\biggr){\tilde a}_1
\int_{\alpha t}^td\tau\int_{\bar\xi\in (a_1,{\tilde a}_1)}
|\Delta_xG(x-\xi,t-\tau)|d\xi.$$

We have proved (3.5), (3.6).

We will estimate now the terms $I_1u$ and $I_4u$
of (3.4).

\vskip 2 mm
{\bf Lemma 8.}
{\it Let function} $u$ {\it satisfy conditions of Theorem} 2ii) {\it and}
$\Delta u$ {\it be represented in} $\Omega_{\sigma}$ {\it by formula} (3.4).
{\it Then terms} $I_1u$ {\it and} $I_4u$ {\it of formula} (3.4) {\it
admit the following (preliminary) estimates for} $\bar x\ge {\tilde a}_1$
{\it and} $t\ge t_0$:
$$\eqalignno{
&|I_1u|\le {4\gamma_0\Gamma^2\cdot {\tilde a}_1^2\over \alpha t}(K^-+K_1),
&(3.9)\cr
&|I_4u|\le {2A_0\Gamma\cdot {\tilde a}_1\over \delta\alpha t}(K^-+K_1),
&(3.10)\cr}$$
{\it where}
$$\eqalign{
&K^-=\int_{\alpha t}^t|\Delta_xG(x-\xi,t-\tau)|_{\scriptstyle
\bar\xi=a_-\atop\scriptstyle \bar\xi< {\tilde a}_1}d\tau,\cr
&K_1=\int_{\alpha t}^t|\Delta_xG(x-\xi,t-\tau)|_{\bar\xi={\tilde a}_1}
d\tau,\cr
&a_-=\bar x\sqrt{t/\tau}-{1\over 2\sqrt{\tau}}-
\sqrt{(t-\tau)/\tau + 1/(4\tau)}.\cr}$$

\vskip 2 mm
{\bf Proof.}
If $t_0$ is large enough and $\tau\ge t_0$ we have using (1.11)
inequalities
$$|u(\xi,\tau)|\le {\Gamma\cdot\bar\xi\over \sqrt{\tau}},\ \
|\v_0(u)|\le 2\gamma_0|u|,\ \ |\Delta_{\xi}\chi(\xi,\tau)|\le
{A_0\over \delta\sqrt{\tau}}.$$
From these relations and from definitions of $I_1u$, $I_4u$ it follows
(using also that $\bar\xi\le {\tilde a}_1\le\bar x$):
$$|I_1u|\le {2\gamma_0\Gamma\cdot {\tilde a}_1\over \sqrt{\alpha t}}I_5u\ \
{\rm and}\ \
|I_4u|\le {A_0\over \delta\sqrt{\alpha t}}I_5u,\eqno(3.11)$$
where
$$I_5u=\int_{\alpha t}^td\tau\int_{\bar\xi\in [a_1,{\tilde a}_1]}
|\Delta_xG(x-\xi,t-\tau)|\cdot |\Delta_{\xi}u(\xi,\tau)|d\xi.\eqno(3.12)$$

The assumption  of Theorem 2ii) implies that
$$\Delta_{\xi}u(\xi,\tau)\ge 0\ \ \forall\tau\ge\tau_0.\eqno(3.13)$$
By Lemma 6 we have also inequalities
$$\eqalign{
&\Delta_xG(x-\xi,t-\tau)<0\ \ {\rm iff}\ \ \bar\xi<\bar x\sqrt{t/\tau},\cr
&\Delta_xG(x-\xi,t-\tau)>0\ \ {\rm iff}\ \ \bar\xi>\bar x\sqrt{t/\tau}.\cr}
\eqno(3.14)$$
From (3.12)-(3.14) we deduce
$$\eqalign{
&I_5u=-\int_{\alpha t}^td\tau\int_{\bar\xi\in [a_1,{\tilde a}_1]}
\Delta_xG\,\Delta_{\xi}ud\xi=\cr
&-\int_{\alpha t}^td\tau\bigl(\int_{\bar\xi\in [a_1,{\tilde a}_1]}
\Delta_x^2G\cdot ud\xi+\Delta_xG\cdot u\big|_{\bar\xi={\tilde a}_1}-
\Delta_xG\cdot u\big|_{\bar\xi=a_1}\bigr).\cr}$$
Using inequality $|u(\xi,\tau)|\le {\Gamma\cdot \bar\xi\over \sqrt{\tau}}$
we obtain
$$|I_5u|\le {\Gamma\cdot {\tilde a}_1\over \sqrt{\alpha t}}
\int_{\alpha t}^td\tau\bigl[\int_{\bar\xi\in [a_1,{\tilde a}_1]}
|\Delta_x^2G|d\xi+|\Delta_xG|_{\bar\xi=a_1}+
|\Delta_xG|_{\bar\xi={\tilde a}_1}\big].\eqno(3.15)$$
From Lemma 6 we have
$$\eqalign{
&\Delta_x^2G(x-\xi,t-\tau)<0,\ \ {\rm iff}\ \ \bar\xi\in (a_-,a_+),\cr
&{\rm where}\ \ a_{\pm}=\bar x\sqrt{t/\tau}-{1\over 2\sqrt{\tau}}\pm
\sqrt{(t-\tau)/\tau +1/(4\tau)}.\cr}\eqno(3.16)$$
If $t_0$ is large enough and ${\tilde a}_1>a_1\sqrt{\alpha}+\sqrt{1-\alpha}$
we have inequality: $a_->a_1$.

Put $\bar\xi_-=\inf\,\{{\tilde a}_1,a_-\}$.

From (3.14), (3.15), (3.16) we deduce
$$\eqalign{
&|I_5u|\le {\Gamma\cdot {\tilde a}_1\over \sqrt{\alpha t}}
\int_{\alpha t}^td\tau\bigl[\int_{\bar\xi\in [a_1,\bar\xi_-]}
\Delta_x^2G\,d\xi-\int_{\bar\xi\in [\bar\xi_-,{\tilde a}_1]}
\Delta_x^2G\,d\xi-\Delta_xG\big|_{\bar\xi=a_1}-
\Delta_xG\big|_{\bar\xi={\tilde a}_1}\bigr]\le\cr
&{\Gamma\cdot {\tilde a}_1\over \sqrt{\alpha t}}
\int_{\alpha t}^td\tau\bigl[
\Delta_xG\big|_{\bar\xi=a_1}-
\Delta_xG\big|_{\bar\xi=\bar\xi_-}+
\Delta_xG\big|_{\bar\xi={\tilde a}_1}-
\Delta_xG\big|_{\bar\xi=\bar\xi_-}-
\Delta_xG\big|_{\bar\xi=a_1}-
\Delta_xG\big|_{\bar\xi={\tilde a}_1}\bigr]\cr
&\le {\Gamma\cdot {\tilde a}_1\over \sqrt{\alpha t}}
\bigl(-2\Delta_xG\big|_{\bar\xi=\bar\xi_-}\bigr)\le
{2\Gamma\cdot {\tilde a}_1\over \sqrt{\alpha t}}(K^-+K_1).\cr}$$
The last  estimate together with
 estimates (3.12) imply (3.9), (3.10).

The following lemma gives more precise estimates for terms
$I_1u$, $I_3u$, $I_4u$.

\vskip 2 mm
{\bf Lemma 9.}
{\it In conditions  and notations of Lemmas} 7,8 {\it we have
estimates}:
$$\eqalignno{
&|I_3u|\le A_2{A_0\Gamma\cdot {\tilde a}_1\over t}
\biggl({\sqrt{1-\alpha}\over \delta^2}+{1\over \delta}+
{{\tilde a}_1\over \delta\sqrt{t}}\biggr),&(3.17)\cr
&|I_1u|\le A_2{\gamma_0\Gamma^2\cdot {\tilde a}_1^2\over  t}
\bigl(1+\ln_+{\sqrt{1-\alpha}\over
{\bar x-{\tilde a}_1}}\bigr),&(3.18)\cr
&|I_4u|\le A_2{A_0\Gamma\cdot {\tilde a}_1\over \delta t}
\bigl(1+\ln_+{\sqrt{1-\alpha}\over
{\bar x-{\tilde a}_1}}\bigr),&(3.19)\cr}$$
{\it where} $A_2$ {\it is absolute constant}, $\alpha$ {\it is sufficiently
close to} 1.

\vskip 2 mm
{\bf Proof.}
In order to prove (3.17), (3.18), (3.19) it is sufficient to prove estimates:
$$\eqalignno{
&K^-\le A\bigl(1+\ln_+{\sqrt{1-\alpha}\over
{\bar x-{\tilde a}_1}}\bigr),&(3.20)\cr
&K_1\le A,&(3.21)\cr
&K^+\le A\sqrt{t}\,\inf\,\{\sqrt{1-\alpha},{1\over {\tilde a}_1}+
{1\over \sqrt{t}}\},&(3.22)\cr}$$
where $K^+$, $K^-$, $K_1$ are integrals from (3.6), (3.9), (3.10).

Let us prove firstly (3.21). Put $\ep=\bar x-{\tilde a}_1$, indicating
that it can be arbitrary small, $y=x-\xi$, $s=t-\tau$. We have
$$\eqalign{
&p=y-s=\bar x\sqrt{t}-{\tilde a}_1\sqrt{\tau}=\cr
&\ep\sqrt{t}+{\tilde a}_1(\sqrt{t}-\sqrt{t-s})=\ep\sqrt{t}+
{{\tilde a}_1s\over 2\theta\sqrt{t}}>0,\cr
&{\rm where}\ \ \theta(s)={{\sqrt{t}+\sqrt{t-s}}\over 2\sqrt{t}},\ \
{{1+\sqrt{\alpha}}\over 2}\le\theta<1.\cr}$$
Since $0\le s\le (1-\alpha)t$ we have
$$\eqalign{
&K_1=\int_{\alpha t}^t|\Delta G|_{\bar\xi={\tilde a}_1}d\tau=-
\int_0^{(1-\alpha)t}\Delta G(p+s,s)ds=K_{10}+K_{11},\cr
&{\rm where}\ \ K_{10}=-\int_{s<p}\Delta G(p+s,s)ds,\ \
K_{11}=-\int_{s>p}\Delta G(p+s,s)ds.\cr}$$
Note that $s<p$ iff $s<\ep\sqrt{t}(1-{{\tilde a}_1\over
2\theta\sqrt{t}})^{-1}$ and ${{\tilde a}_1\over 2\theta\sqrt{t}}<1$.
Hence, inequality $s<p$ implies $s<2\ep\sqrt{t}$, if $t>t_0$ and
inequality $s>p$ implies $s>\ep\sqrt{t}$, if $t>t_0$.

Using Lemma 6ii) we obtain
$$\eqalign{
&K_{10}\le \int_0^{2\ep\sqrt{t}}{1\over \sqrt{s}}e^{-s/4}ds\le\int_0^{\infty}
{1\over \sqrt{s}}e^{-s/4}ds\le A_2\ \ {\rm and}\cr
&K_{11}\le A\int_{\ep\sqrt{t}}^{(1-\alpha)t}s^{-3/2}p\,
e^{-p^2/(4s)}ds\le\cr
&({\rm putting}\ \ s=\eta\cdot t\ \ {\rm and}\ \
p=\sqrt{t}\bigl(\ep+{{\tilde a}_1\eta\over 2\theta}\bigr))\cr
&\le A\int_{\ep/\sqrt{t}}^{1-\alpha}\eta^{-3/2}
\bigl(\ep+{{\tilde a}_1\eta\over 2\theta}\bigr)
\exp\bigl(-\bigl(\ep+{{\tilde a}_1\eta\over 2\theta}\bigr)^2/(4\eta)\bigr)
d\eta\le \cr
&A\bigl(\int_0^1{{\tilde a}_1\over 2\theta\sqrt{\eta}}
\exp\bigl(-{{\tilde a}_1^2\eta\over 16}\bigr)d\eta+
\int_0^1\ep\,\eta^{-3/2}e^{-\ep^2/(4\eta)}d\eta\bigr)\le\cr
&({\rm putting}\ \ \eta=r{\tilde a}_1^{-2}\ \ {\rm or}\ \
\eta=\rho\,\ep^2\ \ {\rm respectively})\cr
&\le A\biggl({1\over 2\theta}\int_0^{\infty}r^{-1/2}e^{-r/16}dr+
\int_0^{\infty}\rho^{-3/2}e^{-1/(4\rho)}d\rho\biggr)\le
A_2.\cr}$$
Inequality (3.21) is proved.

Let us prove now (3.20). Let us find interval of variable $s$ in which
$\bar\xi_-=a_-<{\tilde a}_1$, i.e.
$$\eqalign{
&a_-=\bar x\sqrt{t/\tau}-{1\over 2\sqrt{\tau}}-\sqrt{(t-\tau)/\tau+1/(4\tau)}
<{\tilde a}_1,\cr
&{\rm i.e.}\ \ \bar x\sqrt{t/\tau}-\sqrt{(t-\tau)/\tau}<{\bar a}_1,\ \
{\rm where}\ \ {\bar a}_1={\tilde a}_1\bigl(1+O\bigl({1\over \sqrt{\tau}}
\bigr)\bigr).\cr}$$

Put $\eta={{t-\tau}\over t}={s\over t}$. We obtain
$$\eqalign{
&\bar x-\sqrt{\eta}<{\bar a}_1\sqrt{1-\eta},\ \ {\rm i.e.}\cr
&{\bar x}^2-2{\bar x}\sqrt{\eta}+\eta<{\bar a}_1^2(1-\eta),\ \ {\rm i.e.}\cr
&\bigl(\sqrt{\eta}-{\bar x\over {1+{\bar a}^2}}\bigr)^2<
{{\bar a}_1^2(1+{\bar a}_1^2-{\bar x}^2)\over (1+{\bar a}_1^2)^2},\ \
{\rm i.e.}\cr
&\sqrt{\eta_1}<\sqrt{\eta}<\sqrt{\eta_2},\ \ {\rm where}\cr
&\sqrt{\eta_1}={\bar x\over {1+{\bar a}_1^2}}-
{{\bar a}_1\sqrt{1+{\bar a}_1^2-{\bar x}^2}\over {1+{\bar a}_1^2}};\
\sqrt{\eta_2}={\bar x\over {1+{\bar a}_1^2}}+
{{\bar a}_1\sqrt{1+{\bar a}_1^2-{\bar x}^2}\over {1+{\bar a}_1^2}}.\cr}$$
The interval is not empty if $\bar x\le\sqrt{1+{\bar a}_1^2}$.
In addition we have
$$\eqalign{
&\bar x-{\bar a}_1\sqrt{1+{\bar a}_1^2-{\bar x}^2}\ge
\bar x-{\bar a}_1\bigl(1+{{{\bar a}_1^2-{\bar x}^2}\over 2}\bigr)=\cr
&(\bar x-{\bar a}_1)\bigl(1+{{\bar a}_1(\bar x+{\bar a}_1)\over 2}\bigr)\ge
(\bar x-{\bar a}_1)(1+{\bar a}_1^2).\cr}$$
Hence $\sqrt{\eta_1}>\bar x-{\bar a}_1$.
The condition $\bar\xi_-=a_-$ implies that
$$y=(x-\xi)=(t-\tau)+\sqrt{t-\tau}+O(1)=s+\sqrt{s}+O(1).$$
From Lemmas 5,6 we deduce
$$-\Delta G\big|_{\bar\xi=a_-}\le {\sqrt{e}\over s\sqrt{2\pi}}\bigl(1+
O\bigl({1\over \sqrt{s}}\bigr)\bigr)\le {A\over s}.$$
Hence,
$$K^-\le \int_{\alpha t}^t|\Delta_xG|_{\bar\xi=a_-}d\tau\le
\int_{\eta_1t}^{(1-\alpha)t}{A\over s}ds\le
A_2\ln_+{\sqrt{1-\alpha}\over (\bar x-{\tilde a}_1)},\ \
t\ge t_0.$$

Let us prove (3.22). Using definition of $K^+$ and (3.14) we obtain
$$K^+=-\int_{\alpha t}^td\tau\int_{\bar\xi\in [a_1,{\tilde a}_1]}
\Delta_xG(x-\xi,t-\tau)d\xi\le
\int_{\alpha t}^tG(x-\xi,t-\tau)\big|_{\bar\xi={\tilde a}_1}d\tau.$$
Put (as in the proof of (3.21)) $\ep=\bar x-{\tilde a}_1$,
$y=x-\xi$, $s=t-\tau$, $p=y-s$.

We have
$$\int_{\alpha t}^tG\big|_{\bar\xi={\tilde a}_1}d\tau=
\int_{s<p}G(p+s,s)ds + \int_{s>p}G(p+s,s)ds.$$
Because $s<p$ implies $s<2\ep\sqrt{t}$, $t\ge t_0$ and using Lemma 5i) we
obtain
$$\eqalign{
&\int_{s<p}G(p+s,s)ds\le \int_0^{2\ep\sqrt{t}}{1\over \sqrt{2\pi(p+s)}}
\exp\bigl(-{p^2\over 2(p+s)}\bigr)ds\le\cr
&\int_0^{2\ep\sqrt{t}}{1\over \sqrt{2\pi p}}e^{-p/4}ds\le
\int_0^{2\ep\sqrt{t}}{1\over \sqrt{2\pi s}}e^{-s/4}ds\le A_2.\cr}$$
Because $s>p$ implies $s\in (\ep\sqrt{t},(1-\alpha)t)$ and using Lemma 5i)
we obtain
$$\int_{s>p}G(p+s,s)ds\le \int_{\ep\sqrt{t}}^{(1-\alpha)t}
{1\over \sqrt{2\pi s}}\exp(-{p^2\over 4s})ds\le
{1\over \sqrt{2\pi}}\int_0^{(1-\alpha)t}s^{-1/2}\exp(-{p^2\over 4s})ds.$$
Using $p=\sqrt{t}\bigl(\ep+{{\tilde a}_1s\over 2\theta t}\bigr)$ and
putting $\rho={\tilde a}_1^2{s\over t}$, we obtain further
$$\eqalign{
&\int_{s>p}G(p+s,s)ds\le\cr
&{\sqrt{t}\over \sqrt{2\pi}{\tilde a}_1}\int_0^{{\tilde a}_1^2(1-\alpha)}
{1\over \sqrt{\rho}}e^{-\rho/16}d\rho\le
\sqrt{{t\over 2\pi}}\inf\biggl\{2\sqrt{1-\alpha},
{1\over {\tilde a}_1}\int_0^{\infty}{1\over \sqrt{\rho}}e^{-\rho/16}d\rho
\biggr\}.\cr}$$
Hence, $K^+\le A\sqrt{t}\,\inf\{\sqrt{1-\alpha},{1\over \sqrt{t}}+
{1\over {\tilde a}_1}\}$.

Lemma 9 is proved.

\vskip 2 mm
{\bf Proof of Theorem 2ii).}
From formula (3.4) and estimates (3.5),(3.17),(3.18), (3.19) we deduce
the following inequality under condition that
$\bar x\in ({\tilde a}_1,a_2+\sigma_0\sqrt{t})$, $\sigma_0<\sigma$,
$t\ge {\tilde a}_1^2$ and $\alpha>{{1+\sigma_0}\over {1+\sigma}}$:
$$\eqalign{
&\Delta u(x,t)\le {A_3\Gamma\cdot \bar x\over t}\bigl[{1\over \sqrt{1-\alpha}}
+{\sqrt{1-\alpha}\over \delta^2}+
{1\over \delta}+\cr
&\bigl(\gamma_0\Gamma\cdot {\tilde a}_1 +{1\over \delta}\bigr)
\bigl(1+\ln_+{\sqrt{1-\alpha}\over {\bar x-{\tilde a}_1}}\bigr)\bigr]+\cr
&\gamma_0\Gamma\cdot\bar x\int_{\alpha t}^td\tau\int_{\bar\xi\ge {\tilde a}_1}
|\Delta_xG(x-\xi,t-\tau)|{|\Delta u(\xi,\tau)|\over \sqrt{\tau}}d\xi.\cr}
\eqno(3.23)$$
Put
$$v(t)=t\cdot\max_{\bar x\in ({\tilde a}_1,a_2+\sigma_0\sqrt{t})}
{\Delta u(x,t)\over g(\bar x)},$$
where
$$\eqalign{
&g(\bar x)=B_1+B_2\bigl(1+
\ln_+{\sqrt{1-\alpha}\over {\bar x-{\tilde a}_1}}\bigr),\cr
&B_1={\bar x}\bigl({1\over \sqrt{1-\alpha}}
+{\sqrt{1-\alpha}\over \delta^2}+{1\over \delta}\bigr);\
B_2=\bar x \bigl(\gamma_0\Gamma\cdot {\tilde a}_1 +{1\over \delta}\bigr).\cr}
$$
Then we have $\Delta u(x,t)\le {v(t)\cdot g(\bar x)\over t}$.
From this relation and from (3.23) we obtain
$$v(t)\le A_3\Gamma+{\gamma_0\Gamma\cdot t\over g(\bar x)}\int_{\alpha t}^t
{v(\tau)\over \tau^{3/2}}\int_{\bar\xi>{\tilde a}_1}
|\Delta_xG|\cdot g(\bar\xi)d\xi.$$

By Lemma 6v) we have
$$\int_{\bar\xi>{\tilde a}_1}|\Delta_xG|\cdot g(\bar\xi)d\xi\le
{A_4g(\bar x)\over \sqrt{t-\tau}}
(1+\sqrt{(1-\alpha)/\alpha})(1/\sqrt{\alpha}).$$
From the last two inequalities, putting $\tau=\rho t$, we get
$$v(t)\le A_3\Gamma+A_4\gamma_0\Gamma\int_{\alpha}^1{v(\rho\,t)d\rho\over
\rho^{3/2}\sqrt{1-\rho}}
(1+\sqrt{(1-\alpha)/\alpha})(1/\sqrt{\alpha}).$$
Choose $\alpha_1$ so close to 1 that $\alpha_1>{{1+\sigma_0}\over {1+\sigma}}$
 and
$$(1+\sqrt{(1-\alpha_1)/\alpha_1})(1/\sqrt{\alpha_1})
A_4\gamma_0\Gamma\int_{\alpha_1}^1{d\rho\over
\rho^{3/2}\sqrt{1-\rho}}<1.$$

It means that ${1\over \sqrt{1-\alpha_1}}$ must be of order
$O\bigl({\sqrt{1+\sigma}\over \sqrt{\sigma-\sigma_0}}+\gamma_0\Gamma\bigr)$.
Using Lemma $A_2$ of Appendix we obtain
$$\Delta u\le {v(t)\cdot g(\bar x)\over t}\le {A_5\Gamma\over t}
\biggl(B_1+B_2
\bigl(1+\ln_+{\sqrt{1-\alpha}\over {\bar x-{\tilde a}_1}}\bigr)\biggr),$$
where $\bar x\in ({\tilde a}_1,a_2+\sigma_0\sqrt{t})$, $t\ge t_0\ge
{\tilde a}_1^2$.
Put now $\sqrt{1-\alpha}=\min\,\{\delta,\sqrt{1-\alpha_1}\}$.
Then we obtain
$$\Delta u\le {A_5\Gamma\cdot \bar x\over t}\biggl[{1\over \sqrt{1-\alpha}}+
\bigl(\gamma_0\Gamma\cdot {\tilde a}_1 +{1\over \delta}\bigr)
\bigl(1+\ln_+{\sqrt{1-\alpha}\over {\bar x-{\tilde a}_1}}\bigr)\biggr].$$
Now let $\bar x>a_1$ be fixed and take ${\tilde a}_1={{a_1+\bar x}\over 2}$,
$d={\delta\over 2}$. We obtain
$$\Delta u\le {A_6\Gamma\cdot \bar x\over t}\biggl[
{\sqrt{1+\sigma}\over \sqrt{\sigma-\sigma_0}}+\gamma_0\Gamma+
\bigl(\gamma_0\Gamma\cdot a_1 +{1\over d}\bigr)\biggr].$$
Theorem 2ii) is proved.

\vskip 2 mm
{\bf Sketch of the proof of Theorem} $2ii)^{\prime}$.
\vskip 2 mm
Step 1.
Let function $u$ satisfy equation (1.2) in $\Omega_0$ with
$\v(0)=C=1$ and $\ep=1$. Put $\v_0=\v-C$. We use again the Green-Poisson
type representation formulas for $u$ of type (3.2),(3.4), where
$\chi=\chi_0(\bar x)$, $\bar x={{x-t}\over \sqrt{t}}$, $\chi_0\ :\ \R\to\R$
is a smooth cut-off function such that $0\le\chi_0\le 1$,
$\chi_0\big|_{[{\tilde a}_1,{\tilde a}_2]}\equiv 1$,
$\chi\big|_{(-\infty,a_1)}\equiv 0$, $\chi\big|_{(a_2,\infty)}\equiv 0$,
$0<a_1<{\tilde a}_1<{\tilde a}_2<a_2$, inequalities (3.1) are valid with
$\delta=\min\{{\tilde a}_1-a_1,a_2-{\tilde a}_2\}$. We obtain
representation $(\bar x\in [{\tilde a}_1,{\tilde a}_2],t>\alpha t)$
$$\Delta u=I_0u+I_1u+I_2u+I_3u+I_4u,\eqno(3.4)^{\prime}$$
where
$$\eqalign{
&I_0u=-\int_{\alpha t}^td\tau\int_{\bar\xi\in ({\tilde a}_1,{\tilde a}_2)}
\Delta G\cdot \v_0(u)\cdot \Delta u\,d\xi,\cr
&I_1u=-\int_{\alpha t}^td\tau
\int_{\bar\xi\in [a_1,a_2]\b [{\tilde a}_1,{\tilde a}_2]}
\Delta G\cdot \v_0(u)\cdot \Delta u\,\chi\,d\xi,\cr
&I_2u=\int_{\bar\xi\in [a_1,a_2]}\Delta G(x-\xi,t-\alpha t)
u(\xi,\alpha t)\chi(\xi,\alpha t)\,d\xi,\cr
&I_3u=\int_{\alpha t}^td\tau
\int_{\bar\xi\in [a_1,a_2]\b [{\tilde a}_1,{\tilde a}_2]}
\Delta G\,(u\chi^{\prime}+u\Delta\chi)\,d\xi,\cr
&I_4u=-\int_{\alpha t}^td\tau
\int_{\bar\xi\in [a_1,a_2]\b [{\tilde a}_1,{\tilde a}_2]}
\Delta G\cdot\Delta u\cdot\Delta\chi\,d\xi.\cr}$$

\vskip 2 mm
Step 2. Let $u$ satisfy conditions of Theorem $2ii)^{\prime}$ and
$\Delta u$ be represented  in $\Omega_0$ by formula $(3.4)^{\prime}$,
$\alpha>1/2$.
Using
Lemmas 4,5,6 we obtain Lemma $7^{\prime}$ and $9^{\prime}$:

\vskip 2 mm
{\bf Lemma} $7^{\prime}$.
{\it For} $\bar x\in [a_1,a_2]$ {\it and} $t\ge t_0$
{\it the following estimates are valid}
$$\eqalignno{
&|I_2u(x,t)|\le {A_2\Gamma\over \sqrt{(1-\alpha)}t},&(3.5)^{\prime}\cr
&|I_3u(x,t)|\le {A_2\Gamma\over t}\bigl({1\over \delta^2}+
{{\tilde a}_2\over 2\delta}\bigr)\sqrt{1-\alpha}.
&(3.6)^{\prime}\cr}$$

\vskip 2 mm
{\bf Lemma} $9^{\prime}$.
{\it For} $\bar x\in [{\tilde a}_1,{\tilde a}_2]$  {\it and} $t\ge t_0$
{\it the following estimates are valid}
$$\eqalignno{
&|I_1u(x,t)|\le A_2{\gamma_0\Gamma^2\over  t}
\bigl(1+\ln_+{\sqrt{1-\alpha}\over {\bar x-{\tilde a}_1}}+
\ln_+{\sqrt{1-\alpha}\over {{\tilde a}_2-\bar x}}
\bigr),&(3.18)^{\prime}\cr
&|I_4u(x,t)|\le A_2{A_0\Gamma\over \delta t}
\bigl(1+\ln_+{\sqrt{1-\alpha}\over {\bar x-{\tilde a}_1}}+
\ln_+{\sqrt{1-\alpha}\over {{\tilde a}_2-\bar x}}
\bigr).&(3.19)^{\prime}\cr}$$

\vskip 2 mm
Step 3. From formula $(3.4)^{\prime}$ and estimates $(3.5)^{\prime}$,
$(3.6)^{\prime}$, $(3.18)^{\prime}$, $(3.19)^{\prime}$  we deduce the
following inequality ($\bar x\in [{\tilde a}_1,{\tilde a}_2]$)
$$\eqalign{
&\Delta u\le {A_3\Gamma\over t}\bigg[{1\over \sqrt{1-\alpha}}+
{\sqrt{1-\alpha}\over \delta^2}+ {{\tilde a}_2\sqrt{1-\alpha}\over 2\delta}+
\cr
&\bigl(\gamma_0\Gamma+{1\over \delta}\bigr)
\bigl(1+\ln_+{\sqrt{1-\alpha}\over {\bar x-{\tilde a}_1}}+
\ln_+{\sqrt{1-\alpha}\over {{\tilde a}_2-\bar x}}
\bigr)\biggr]-\cr
&\int_{\alpha t}^td\tau\int_{\bar\xi\in ({\tilde a}_1,{\tilde a}_2)}
\Delta G\,\v_0(u)\,\Delta u\,d\xi.\cr}\eqno(3.23)^{\prime}$$

By assumption of Theorem $2ii)^{\prime}$ we have
$\Delta_{\xi}(\xi,\tau)\ge 0$. If in
 assumptions of Theorem $2ii)^{\prime}$ we have additional positivity
conditions $\v^{\prime}(0)\ge 0$ and $u\ge 0$
then we can replace the integral term in $(3.23)^{\prime}$ by the
following bigger one
$$-\gamma_0\Gamma\int_{\alpha t}^t{d\tau\over \sqrt{\tau}}
\int_{\bar\xi\ :\ \Delta G<0}\Delta_xG\cdot\Delta u\,d\xi.$$
Following further the proof of Theorem 2ii)  and applying again Lemma 6v)
we obtain the statement of Theorem $2ii)^{\prime}$ with constant
$B=B_0(a_2+{1\over d}+{\gamma_0\Gamma\over C})$.

Without additional positivity conditions the statement of Theorem
$2ii)^{\prime}$ is also valid but for the proof of it more hard version of
Lemma 6v) is needed
where the weight $\bigl(1+\ln_+{1\over {\bar\xi-{\tilde a}_1}}\bigr)$ is
replaced by
$\bigl(1+\ln_+{1\over {{\tilde a}_2-\bar\xi}}\bigr)$.

\vskip 2 mm
{\bf Lemma} 6v)$^{\prime}$.
Let $0<\bar x<{\tilde a}_2$. Then
$$\int_{\bar\xi<{\tilde a}_2}|\Delta G(x-\xi,t-\tau)|\biggl(1+
\ln_+{1\over {{\tilde a}_2-\bar\xi}}\biggr)d\xi\le
{A_1^{\prime}\over \sqrt{t-\tau}}\biggl(1+\ln_+{\tilde a}_2+
\ln_+{1\over {{\tilde a}_2-\bar x}}\biggr).$$

\vskip 4 mm
{\bf Appendix.} Integral inequalities.

\vskip 2 mm
{\bf Lemma} $A_1$.
{\it Let} $0\le\psi(x)=O\bigl({1\over x}\bigr)$, $x\ge 0$, {\it and}
$\int_0^{\infty}\psi(x)dx<\infty$. {\it Then}
$$\int_0^a\psi(x)\ln_+{b\over {a-x}}dx\le A_{\psi}\bigl(1+\ln_+{b\over a}
\bigr).$$

\vskip 2 mm
{\bf Proof.}
Let $a<b$. Then
$$\eqalign{
&\int_0^a\psi(x)\ln_+{b\over {a-x}}dx=\int_0^{a/2}\psi(x)\ln_+{b\over {a-x}}dx
+ \int_{a/2}^a\psi(x)\ln_+{b\over {a-x}}dx\le\cr
&\ln_+{2b\over a}\int_0^{\infty}\psi(x)dx+\max_{x> a/2}\psi(x)\int_0^a
\ln_+{b\over {a-x}}dx=\cr
&A_{\psi}\bigl({1\over 2}\ln_+{2b\over a}+\ln_+{b\over a}+1\bigr)\le
A_{\psi}\bigl(\ln_+{b\over a}+1\bigr).\cr}$$
Let $a>b$. Then
$$\eqalign{
&\int_0^a\psi(x)\ln_+{b\over {a-x}}dx=
\int_{a-b}^a\psi(x)\ln_+{b\over {a-x}}dx=\cr
&\int_{a-b}^{a-b/2}\psi(x)\ln_+{b\over {a-x}}dx+
\int_{a-b/2}^a\psi(x)\ln_+{b\over {a-x}}dx\le \cr
&\ln_+2\int_0^{\infty}\psi(x)dx+\max_{x> a/2}\psi(x)\int_0^{b/2}
\ln_+{b\over x}dx\le A_{\psi}.\cr}$$

\vskip 2 mm
{\bf Lemma} $A_2$.
{\it Let} $v(t)$ {\it be a continuous function satisfying the inequality}
$$v(t)\le A+\int_{\alpha}^1h(\rho)v(\rho\,t)d\rho,\ \ t\ge t_0,$$
{\it where}
$$0<\int_{\alpha}^1h(\rho)d\rho<1,\ \ h\ge 0,\ \ \alpha\in (0,1).$$
{\it Then}\ \ $\exists\ m>0,\ M>0$ {\it such that}
$v(t)\le A_1+Mt^{-m}$, $t\ge t_0$, {\it where}
$$A_1=A\bigl(1-\int_{\alpha}^1h(\rho)d\rho\bigr)^{-1}.$$

\vskip 2 mm
{\bf Proof.}
Find $A_1\in\R$ such that $v_1(t)=A_1$ satisfies the equation
$$v_1(t)=A+\int_{\alpha}^1h(\rho)v_1(\rho t)d\rho.$$
We get
$$A_1=A\bigl(1-\int_{\alpha}^1h(\rho)d\rho\bigr)^{-1}.$$
Let us find $m>0$ such that $v_0(t)=1/t^m$ satisfies the equation
$$v_0(t)=\int_{\alpha}^1h(\rho)v_0(\rho t)d\rho.$$
This holds iff $\int_{\alpha}^1{h(\rho)\over \rho^m}d\rho=1$.
Since $I(m)=\int_{\alpha}^1{h(\rho)\over \rho^m}d\rho$ is a continuous
function of $m$, $I(m)\to +\infty$ as $m\to +\infty$, $I(0)<1$, then
there exists $m$ such that $I(m)=1$.

Choose $M$ large enough such that
$$V(t)=v(t)-v_1(t)-Mv_0(t)<0$$
for $t_0<t\le t_0/\alpha =t_1$.
We claim that $V(t)<0\ \ \forall\ t\ge t_0$.

Indeed, let $t^*=\sup\,\{t\ge t_0\ :\ V(t)<0\}$. By the choice of $M$ and
continuity of $V$ we have $t^*>t_1$.

If $t^*$ is finite then
$$V(t^*)\le \int_{\alpha}^1h(\rho)V(\rho t^*)d\rho < 0.$$
Since $V$ is continuous, $V<0$ holds in a neighborhood of $t^*$, but this
contradicts to the definition of $t^*$.

\vfill\eject

\centerline{\bf References.}
\vskip 4 mm
\item{[ 1]} Bateman H., Some recent researches on the motion of fluids,
Monthly Weather Review, {\bf 43}, 1915, 163-170
\item{[ 2]} Belenky V., Diagram of growth of a monotonic function and a
problem of their reconstruction by the Diagram, Preprint, Central
Economics and Mathematical Institute, Academy of Sciences of the USSR,
Moscow, 1990, 1-44 (in Russian)
\item{[ 3]} Burgers J.M., Application of a model system to illustrate
some points of the statistical theory of free turbulence, Proc. Acad. Sci.
Amsterdam, {\bf 43}, 1940, 2-12
\item{[ 4]} Gelfand I.M., Some problems in the theory of quasilinear
equations, Usp. Mat. Nauk {\bf 14}, 1959, 87-158 (in Russian); Amer. Math.
Soc. Translations, {\bf 33}, 1963
\item{[ 5]} Henkin G.M., Polterovich V.M., Schumpeterian dynamics as a
nonlinear wave theory, J.Math. Econ. {\bf 20}, 1991, 551-590
\item{[ 6]} Henkin G.M., Polterovich V.M., A difference-differential
analogue of the Burgers equation and some models of economic
development, Discrete Contin. Dyn. Syst., {\bf 5}, 1999, 697-728
\item{[ 7]} Henkin G.M., Shananin A.A., Asymptotic behavior of solutions
of the Cauchy problem for Burgers type equations, Preprint, 2004
\item{[ 8]} Hopf E., The partial differential equation
$u_t+uu_x=\mu u_{xx}$, Comm. on Pure and Appl. Math., {\bf 3}, 1950,
201-230
\item{[ 9]} Iljin A.M., Olejnik O.A., Asymptotic long-time behavior
of the Cauchy problem for some quasilinear equation, Mat. Sbornik, {\bf 51},
1960, 191-216 (in Russian)
\item{[10]} Landau L.D., Lifchitz E.M., M\'ecanique des fluides, 2 \`eme
\'edition, MIR, Moscou, 1989
\item{[11]} Lui T.-P., Matsumura A., Nishihara K., Behaviors of solutions
for the Burgers equation with boundary corresponding to rarefaction
waves, SIAM J.Math.Anal., {\bf 29}, 1998, 293-308
\item{[12]} Oleinik O.A., Uniqueness and stability of the generalized
solution of the Cauchy problem for a quasilinear equation, Usp. Mat. Nauk
{\bf 14}, 1959, 165-170 (in Russian); Amer. Math. Soc. Translations {\bf 33},
1963, 285-290
\item{[13]} Serre D., $L^1$-stability of nonlinear waves in scalar
conservation laws, Handbook of Differential Equations, C.Dafermos,
E.Feireisl eds, Elsevier, 2004
\item{[14]} Weinberger H.F., Long-time behavior for a regularized scalar
conservation law in the absence of genuine nonlinearity, Ann. Inst. H.
Poincare, Analyse Nonlineaire, 1990, 407-425

\end